\documentclass[11pt]{article}
\usepackage{latexsym,amssymb,amsmath,geometry,cite,enumerate,float,verbatim}
\newtheorem{theorem}{Theorem}[section]
\newtheorem{lemma}[theorem]{Lemma}
\newtheorem{observation}[theorem]{Observation}

\usepackage{lineno}
\usepackage{setspace}

\begin{document}

%\linenumbers
\onehalfspace

\title{Sandwiches Missing Two Ingredients of Order Four}

\author{Jos\'{e} D. Alvarado$^1$\and 
Simone Dantas$^1$\and
Dieter Rautenbach$^2$}

\date{}

\maketitle

\begin{center}
{\small 
$^1$ Instituto de Matem\'{a}tica e Estat\'{i}stica, Universidade Federal Fluminense, Niter\'{o}i, Brazil\\
\texttt{josealvarado.mat17@gmail.com, sdantas@im.uff.br}\\[3mm]
$^2$ Institute of Optimization and Operations Research, Ulm University, Ulm, Germany\\
\texttt{dieter.rautenbach@uni-ulm.de}
}
\end{center}

\begin{abstract}
For a set ${\cal F}$ of graphs, 
an instance of the ${\cal F}$-{\sc free Sandwich Problem} is a pair $(G_1,G_2)$
consisting of two graphs $G_1$ and $G_2$ with the same vertex set such that $G_1$ is a subgraph of $G_2$,
and the task is to determine an ${\cal F}$-free graph $G$ containing $G_1$ and contained in $G_2$,
or to decide that such a graph does not exist.
Initially motivated by the graph sandwich problem for trivially perfect graphs,
which are the $\{ P_4,C_4\}$-free graphs,
we study the complexity of the ${\cal F}$-{\sc free Sandwich Problem} for sets ${\cal F}$ containing two non-isomorphic graphs of order four.
We show that
if ${\cal F}$ is one of the sets 
$\left\{ {\rm diamond},K_4\right\}$,
$\left\{ {\rm diamond},C_4\right\}$,
$\left\{ {\rm diamond},{\rm paw}\right\}$,
$\left\{ K_4,\overline{K_4}\right\}$,
$\left\{ P_4,C_4\right\}$,
$\left\{ P_4,\overline{\rm claw}\right\}$,
$\left\{ P_4,\overline{\rm paw}\right\}$,
$\left\{ P_4,\overline{\rm diamond}\right\}$,
$\left\{ {\rm paw},C_4\right\}$,
$\left\{ {\rm paw},{\rm claw}\right\}$,
$\left\{ {\rm paw},\overline{{\rm claw}}\right\}$,
$\left\{ {\rm paw},\overline{\rm paw}\right\}$,
$\left\{ C_4,\overline{C_4}\right\}$,
$\left\{ {\rm claw},\overline{{\rm claw}}\right\}$, and
$\left\{ {\rm claw},\overline{C_4}\right\}$,
then the ${\cal F}$-{\sc free Sandwich Problem} can be solved in polynomial time, and,
if ${\cal F}$ is one of the sets
$\left\{ C_4,K_4\right\}$,
$\left\{ {\rm paw},K_4\right\}$,
$\left\{ {\rm paw},\overline{K_4}\right\}$,
$\left\{ {\rm paw},\overline{C_4}\right\}$,
$\left\{ {\rm diamond},\overline{C_4}\right\}$,
$\left\{ {\rm paw},\overline{\rm diamond}\right\}$, and
$\left\{ {\rm diamond},\overline{\rm diamond}\right\}$,
then the decision version of the ${\cal F}$-{\sc free Sandwich Problem} is NP-complete.
\end{abstract}

{\small 

\begin{tabular}{lp{13cm}}
{\bf Keywords:} Graph sandwich problem; forbidden induced subgraph
\end{tabular}
}

\pagebreak

\section{Introduction}

Graph sandwich problems \cite{gokash} are a natural generalization of recognition problems, 
and have received considerable attention \cite{ceevfikl,dafisite,dafimate,dafigoklma,daklmemo,go,tedafi,tedafi2}.
It is not unusual that graph classes for which the recognition is very easy lead to challenging graph sandwich problems,
which are either intractable or require interesting structural and algorithmic arguments for their solution.
In such a situation, the graph sandwich problem motivates a detailed analysis of the corresponding graph class leading to insights
that were probably not needed for some efficient ad-hoc recognition algorithm 
but are essential for the solution of the sandwich problem.

Good examples for this effect are graph classes defined by a finite set ${\cal F}$ of forbidden induced subgraphs.
In \cite{dafisite} Dantas, de Figueiredo, da Silva, and Teixeira 
initiated the study of graph sandwich problems for ${\cal F}$-free graphs, where ${\cal F}$ contains a single graph.
In \cite{dafimate} Dantas, de Figueiredo, Maffray, and Teixeira provided further results along this line,
and, in particular, settled the complexity status of the graph sandwich problem for $\{ F\}$-free graphs
for every graph $F$ of order four.
Considering forbidden induced subgraph of order four is rather natural, 
because many well known graph classes \cite{brlesp} are defined by one or more such graphs,
and various aspects of these classes have been studied \cite{brenlelo,dapa,kolo}.

Originally motivated by the graph sandwich problem for trivially perfect graphs,
which are the $\{ P_4,C_4\}$-free graphs,
and following a suggestion by Golumbic, 
we initiate the study of the graph sandwich problem for ${\cal F}$-free graphs,
where ${\cal F}$ is a set of two non-isomorphic graphs of order four.
In order to obtain our results, we rely on known results \cite{brma,br,ol,mapr,mape} for some cases,
and develop new arguments for other cases.

\medskip

\noindent Before we proceed to our results, we recall some relevant definitions.
We consider finite, simple, and undirected graphs, and use standard terminology and notation.
For a graph property $\Pi$, that is, $\Pi$ is a set of graphs, the corresponding graph sandwich problem is the following.

\medskip

\noindent $\Pi$-{\sc Sandwich Problem}\\
\begin{tabular}{lp{14cm}}
Instance: & A pair $(G_1,G_2)$ of two graphs such that $G_1$ and $G_2$ have the same vertex set, and $G_1$ is a subgraph of $G_2$.\\
Task: & Determine a graph $G$ with $G_1\subseteq G\subseteq G_2$ and $G\in \Pi$, or conclude that no such graph exists.
\end{tabular}

\medskip

\noindent Let ${\cal F}$ be a set of graphs.
A graph $G$ is ${\cal F}$-free if no induced subgraph of $G$ is in ${\cal F}$.
Let $\overline{\cal F}$ be $\left\{ \overline{F}:F\in {\cal F}\right\}$, 
where $\overline{F}$ is the complement of a graph $F$.
For two graphs $G_1$ and $G_2$ such that 
$G_1$ and $G_2$ have the same vertex set, and $G_1$ is a subgraph of $G_2$, 
let ${\cal SW}_{\cal F}(G_1,G_2)$
be the set of ${\cal F}$-free graphs $G$ with $G_1\subseteq G\subseteq G_2$.

Here is the type of problem we consider.

\medskip

\noindent ${\cal F}$-{\sc free Sandwich Problem}\\
\begin{tabular}{lp{14cm}}
Instance: & A pair $(G_1,G_2)$ of two graphs such that $G_1$ and $G_2$ have the same vertex set, and $G_1$ is a subgraph of $G_2$.\\
Task: & Determine a graph $G$ in ${\cal SW}_{\cal F}(G_1,G_2)$, or conclude that this set is empty.
\end{tabular}

\medskip

\noindent The ${\cal F}$-{\sc free Sandwich Decision Problem} has the same input as the ${\cal F}$-{\sc free Sandwich Problem}
but the task is merely to decide whether ${\cal SW}_{\cal F}(G_1,G_2)$ is non-empty.
It is easy to see that the ${\cal F}$-{\sc free Sandwich Problem} can be solved in polynomial time
if and only if the ${\cal F}$-{\sc free Sandwich Decision Problem} can.
In fact, if ${\cal SW}_{\cal F}(G_1,G_2)$ is non-empty, 
then iteratively applying an efficient algorithm for the ${\cal F}$-{\sc free Sandwich Decision Problem},
one can determine in polynomial time an edge-minimal graph $G$ with $G_1\subseteq G\subseteq G_2$
such that ${\cal SW}_{\cal F}(G_1,G)$ is still non-empty,
and this graph $G$ actually lies in ${\cal SW}_{\cal F}(G_1,G_2)$.
We collect some simple observations.

\begin{observation}\label{observation1}
Let ${\cal F}$ be a set of graphs, and let $(G_1,G_2)$ be an instance of the ${\cal F}$-{\sc free Sandwich Problem}.
\begin{enumerate}[(i)]
\item ${\cal SW}_{\overline{\cal F}}\left(\overline{G_2},\overline{G_1}\right)=\overline{{\cal SW}_{\cal F}(G_1,G_2)}$.
\item If all graphs in ${\cal F}$ are connected, and ${\cal SW}_{\cal F}(G_1,G_2)$ is non-empty,
then there is some graph $G$ in ${\cal SW}_{\cal F}(G_1,G_2)$ 
such that the vertex sets of the components of $G_1$ are the same as the vertex sets of the components of $G$.
\item If no graph in ${\cal F}$ has a universal vertex, and $u$ is a universal vertex in $G_2$,
then ${\cal SW}_{\cal F}(G_1,G_2)$ is non-empty 
if and only if ${\cal SW}_{\cal F}(G_1-u,G_2-u)$ is non-empty.
\item If every graph $F$ in ${\cal F}$
has a unique ${\cal F}$-free supergraph $F^*$ with $V(F)=V(F^*)$,
then the ${\cal F}$-{\sc free Sandwich Problem} can be solved in polynomial time.
\end{enumerate}
\end{observation}
{\it Proof:} (i) This follows immediately from the definition.

\medskip

\noindent (ii) Since the vertex set of each component of a graph $G$ in ${\cal SW}_{\cal F}(G_1,G_2)$
is the union of vertex sets of components of $G_1$, 
and all edges of $G$ between components of $G_1$ belong to $G_2$, 
removing from $G$ all such edges
yields another graph in ${\cal SW}_{\cal F}(G_1,G_2)$ that has the desired property.

\medskip

\noindent (iii) If $G\in {\cal SW}_{\cal F}(G_1,G_2)$, then $G-u\in {\cal SW}_{\cal F}(G_1-u,G_2-u)$,
which implies the necessity.
By the assumption on ${\cal F}$, 
adding a universal vertex to an ${\cal F}$-free graph yields an ${\cal F}$-free graph,
which implies the sufficiency.

\medskip

\noindent (iv) Starting with $G_1$, 
and iteratively adding the uniquely determined sets of edges to every induced subgraph from ${\cal F}$ using edges of $G_2$
yields a graph in ${\cal SW}_{\cal F}(G_1,G_2)$.
If, at some point, the graph $G_2$ does not contain the necessary edges, then ${\cal SW}_{\cal F}(G_1,G_2)$ is empty. $\Box$

\medskip

\noindent As said above, our goal it to study the complexity of the ${\cal F}$-{\sc free Sandwich Problem}
for sets ${\cal F}$ containing two non-isomorphic graphs of order four.
Figure \ref{fig0} illustrates all such graphs together with the names we are using.
By Observation \ref{observation1}(i), 
it suffices to consider the sets ${\cal F}$ up to complementation.
Note that $P_4$ is the only self-complementary graph of order four.
Hence, up to complementation, there are $5$ sets ${\cal F}$ that contain $P_4$.
There are $10$ sets ${\cal F}$ containing two non-isomorphic graphs with less than four edges,
and, up to complementation, there are $15$ sets ${\cal F}$ 
containing one graph with less than four edges
and one graph with more than four edges.
Altogether, the $30$ choices for ${\cal F}$ illustrated in Figure \ref{fig1} 
represent all sets of two non-isomorphic graphs of order four up to complementation.

\begin{figure}
\begin{center}
%TeXCAD Picture [3.pic]. Options:
%\grade{\on}
%\emlines{\off}
%\epic{\off}
%\beziermacro{\on}
%\reduce{\on}
%\snapping{\on}
%\pvinsert{% Your \input, \def, etc. here}
%\quality{8.000}
%\graddiff{0.005}
%\snapasp{1}
%\zoom{8.0001}
\unitlength 0.8mm % = 2.845pt
\linethickness{0.4pt}
\ifx\plotpoint\undefined\newsavebox{\plotpoint}\fi % GNUPLOT compatibility
\begin{picture}(161,41)(0,0)
\put(30,30){\circle*{2}}
\put(0,30){\circle*{2}}
\put(60,30){\circle*{2}}
\put(90,30){\circle*{2}}
\put(120,30){\circle*{2}}
\put(150,30){\circle*{2}}
\put(15,5){\circle*{2}}
\put(45,5){\circle*{2}}
\put(75,5){\circle*{2}}
\put(105,5){\circle*{2}}
\put(135,5){\circle*{2}}
\put(40,30){\circle*{2}}
\put(10,30){\circle*{2}}
\put(70,30){\circle*{2}}
\put(100,30){\circle*{2}}
\put(130,30){\circle*{2}}
\put(160,30){\circle*{2}}
\put(25,5){\circle*{2}}
\put(55,5){\circle*{2}}
\put(85,5){\circle*{2}}
\put(115,5){\circle*{2}}
\put(145,5){\circle*{2}}
\put(40,40){\circle*{2}}
\put(10,40){\circle*{2}}
\put(70,40){\circle*{2}}
\put(100,40){\circle*{2}}
\put(130,40){\circle*{2}}
\put(160,40){\circle*{2}}
\put(25,15){\circle*{2}}
\put(55,15){\circle*{2}}
\put(85,15){\circle*{2}}
\put(115,15){\circle*{2}}
\put(145,15){\circle*{2}}
\put(30,40){\circle*{2}}
\put(0,40){\circle*{2}}
\put(60,40){\circle*{2}}
\put(90,40){\circle*{2}}
\put(120,40){\circle*{2}}
\put(150,40){\circle*{2}}
\put(15,15){\circle*{2}}
\put(45,15){\circle*{2}}
\put(75,15){\circle*{2}}
\put(105,15){\circle*{2}}
\put(135,15){\circle*{2}}
\put(40,40){\line(0,1){0}}
\put(10,40){\line(0,1){0}}
\put(70,40){\line(0,1){0}}
\put(100,40){\line(0,1){0}}
\put(130,40){\line(0,1){0}}
\put(160,40){\line(0,1){0}}
\put(25,15){\line(0,1){0}}
\put(55,15){\line(0,1){0}}
\put(85,15){\line(0,1){0}}
\put(115,15){\line(0,1){0}}
\put(145,15){\line(0,1){0}}
\put(30,40){\line(1,0){10}}
\put(150,40){\line(1,0){10}}
\put(45,15){\line(1,0){10}}
\put(75,15){\line(1,0){10}}
\put(105,15){\line(1,0){10}}
\put(135,15){\line(1,0){10}}
\put(60,40){\line(0,-1){10}}
\put(60,30){\line(1,1){10}}
\put(90,40){\line(0,-1){11}}
\put(100,40){\line(0,-1){10}}
\put(130,40){\line(-1,0){11}}
\put(120,40){\line(0,-1){10}}
\put(120,30){\line(1,1){10}}
\put(150,30){\line(0,1){10}}
\put(160,30){\line(0,1){10}}
\put(25,5){\line(-1,0){10}}
\put(15,5){\line(1,1){10}}
\put(15,15){\line(0,-1){10}}
\put(45,15){\line(0,-1){10}}
\put(45,5){\line(1,0){10}}
\put(55,15){\line(-1,-1){10}}
\put(75,15){\line(0,-1){10}}
\put(75,5){\line(1,0){10}}
\put(85,5){\line(0,1){10}}
\put(105,15){\line(0,-1){10}}
\put(135,15){\line(0,-1){10}}
\put(105,5){\line(1,0){10}}
\put(135,5){\line(1,0){10}}
\put(115,5){\line(0,1){10}}
\put(145,5){\line(0,1){10}}
\put(115,15){\line(-1,-1){10}}
\put(145,15){\line(-1,-1){10}}
\put(135,15){\line(1,-1){10}}
\put(20,0){\makebox(0,0)[cc]{claw}}
\put(50,0){\makebox(0,0)[cc]{paw}}
\put(80,0){\makebox(0,0)[cc]{$C_4$}}
\put(110,0){\makebox(0,0)[cc]{diamond}}
\put(35,25){\makebox(0,0)[cc]{$\overline{\rm diamond}$}}
\put(5,25){\makebox(0,0)[cc]{$\overline{K_4}$}}
\put(65,25){\makebox(0,0)[cc]{$\overline{\rm paw}$}}
\put(95,25){\makebox(0,0)[cc]{$\overline{C_4}$}}
\put(125,25){\makebox(0,0)[cc]{$\overline{\rm claw}$}}
\put(155,25){\makebox(0,0)[cc]{$P_4$}}
\put(140,0){\makebox(0,0)[cc]{$K_4$}}
\end{picture}
\end{center}
\caption{All graphs of order four.}\label{fig0}
\end{figure}
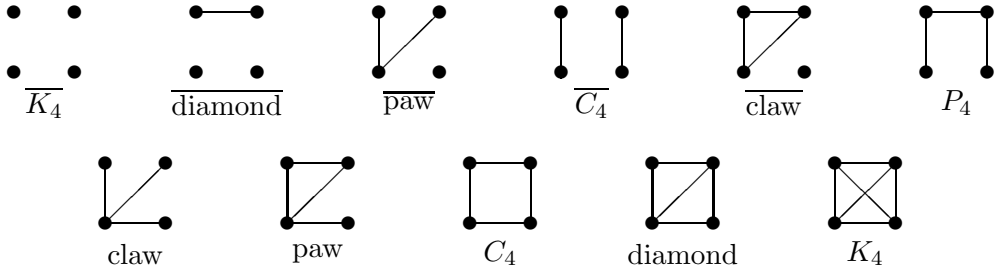

\begin{figure}[H]
\begin{center}
%TeXCAD Picture [2.pic]. Options:
%\grade{\on}
%\emlines{\off}
%\epic{\off}
%\beziermacro{\on}
%\reduce{\on}
%\snapping{\on}
%\pvinsert{% Your \input, \def, etc. here}
%\quality{8.000}
%\graddiff{0.005}
%\snapasp{1}
%\zoom{5.6569}
\unitlength 0.87mm % = 2.845pt
\linethickness{0.4pt}
\ifx\plotpoint\undefined\newsavebox{\plotpoint}\fi % GNUPLOT compatibility
\begin{picture}(196,193)(0,0)
\put(4,4){\circle*{2}}
\put(22,4){\circle*{2}}
\put(22,43){\circle*{2}}
\put(42,4){\circle*{2}}
\put(42,43){\circle*{2}}
\put(42,82){\circle*{2}}
\put(62,4){\circle*{2}}
\put(62,43){\circle*{2}}
\put(62,82){\circle*{2}}
\put(62,121){\circle*{2}}
\put(82,4){\circle*{2}}
\put(82,43){\circle*{2}}
\put(82,82){\circle*{2}}
\put(82,121){\circle*{2}}
\put(82,160){\circle*{2}}
\put(102,4){\circle*{2}}
\put(102,43){\circle*{2}}
\put(102,82){\circle*{2}}
\put(102,121){\circle*{2}}
\put(102,160){\circle*{2}}
\put(122,4){\circle*{2}}
\put(122,43){\circle*{2}}
\put(122,82){\circle*{2}}
\put(122,121){\circle*{2}}
\put(142,4){\circle*{2}}
\put(142,43){\circle*{2}}
\put(142,82){\circle*{2}}
\put(162,4){\circle*{2}}
\put(162,43){\circle*{2}}
\put(182,4){\circle*{2}}
\put(4,20){\circle*{2}}
\put(22,20){\circle*{2}}
\put(22,59){\circle*{2}}
\put(42,20){\circle*{2}}
\put(42,59){\circle*{2}}
\put(42,98){\circle*{2}}
\put(62,20){\circle*{2}}
\put(62,59){\circle*{2}}
\put(62,98){\circle*{2}}
\put(62,137){\circle*{2}}
\put(82,20){\circle*{2}}
\put(82,59){\circle*{2}}
\put(82,98){\circle*{2}}
\put(82,137){\circle*{2}}
\put(82,176){\circle*{2}}
\put(102,20){\circle*{2}}
\put(102,59){\circle*{2}}
\put(102,98){\circle*{2}}
\put(102,137){\circle*{2}}
\put(102,176){\circle*{2}}
\put(122,20){\circle*{2}}
\put(122,59){\circle*{2}}
\put(122,98){\circle*{2}}
\put(122,137){\circle*{2}}
\put(142,20){\circle*{2}}
\put(142,59){\circle*{2}}
\put(142,98){\circle*{2}}
\put(162,20){\circle*{2}}
\put(162,59){\circle*{2}}
\put(182,20){\circle*{2}}
\put(14,4){\circle*{2}}
\put(32,4){\circle*{2}}
\put(32,43){\circle*{2}}
\put(52,4){\circle*{2}}
\put(52,43){\circle*{2}}
\put(52,82){\circle*{2}}
\put(72,4){\circle*{2}}
\put(72,43){\circle*{2}}
\put(72,82){\circle*{2}}
\put(72,121){\circle*{2}}
\put(92,4){\circle*{2}}
\put(92,43){\circle*{2}}
\put(92,82){\circle*{2}}
\put(92,121){\circle*{2}}
\put(92,160){\circle*{2}}
\put(112,4){\circle*{2}}
\put(112,43){\circle*{2}}
\put(112,82){\circle*{2}}
\put(112,121){\circle*{2}}
\put(112,160){\circle*{2}}
\put(132,4){\circle*{2}}
\put(132,43){\circle*{2}}
\put(132,82){\circle*{2}}
\put(132,121){\circle*{2}}
\put(152,4){\circle*{2}}
\put(152,43){\circle*{2}}
\put(152,82){\circle*{2}}
\put(172,4){\circle*{2}}
\put(172,43){\circle*{2}}
\put(192,4){\circle*{2}}
\put(14,20){\circle*{2}}
\put(32,20){\circle*{2}}
\put(32,59){\circle*{2}}
\put(52,20){\circle*{2}}
\put(52,59){\circle*{2}}
\put(52,98){\circle*{2}}
\put(72,20){\circle*{2}}
\put(72,59){\circle*{2}}
\put(72,98){\circle*{2}}
\put(72,137){\circle*{2}}
\put(92,20){\circle*{2}}
\put(92,59){\circle*{2}}
\put(92,98){\circle*{2}}
\put(92,137){\circle*{2}}
\put(92,176){\circle*{2}}
\put(112,20){\circle*{2}}
\put(112,59){\circle*{2}}
\put(112,98){\circle*{2}}
\put(112,137){\circle*{2}}
\put(112,176){\circle*{2}}
\put(132,20){\circle*{2}}
\put(132,59){\circle*{2}}
\put(132,98){\circle*{2}}
\put(132,137){\circle*{2}}
\put(152,20){\circle*{2}}
\put(152,59){\circle*{2}}
\put(152,98){\circle*{2}}
\put(172,20){\circle*{2}}
\put(172,59){\circle*{2}}
\put(192,20){\circle*{2}}
\put(14,14){\circle*{2}}
\put(32,14){\circle*{2}}
\put(32,53){\circle*{2}}
\put(52,14){\circle*{2}}
\put(52,53){\circle*{2}}
\put(52,92){\circle*{2}}
\put(72,14){\circle*{2}}
\put(72,53){\circle*{2}}
\put(72,92){\circle*{2}}
\put(72,131){\circle*{2}}
\put(92,14){\circle*{2}}
\put(92,53){\circle*{2}}
\put(92,92){\circle*{2}}
\put(92,131){\circle*{2}}
\put(92,170){\circle*{2}}
\put(112,14){\circle*{2}}
\put(112,53){\circle*{2}}
\put(112,92){\circle*{2}}
\put(112,131){\circle*{2}}
\put(112,170){\circle*{2}}
\put(132,14){\circle*{2}}
\put(132,53){\circle*{2}}
\put(132,92){\circle*{2}}
\put(132,131){\circle*{2}}
\put(152,14){\circle*{2}}
\put(152,53){\circle*{2}}
\put(152,92){\circle*{2}}
\put(172,14){\circle*{2}}
\put(172,53){\circle*{2}}
\put(192,14){\circle*{2}}
\put(14,30){\circle*{2}}
\put(32,30){\circle*{2}}
\put(32,69){\circle*{2}}
\put(52,30){\circle*{2}}
\put(52,69){\circle*{2}}
\put(52,108){\circle*{2}}
\put(72,30){\circle*{2}}
\put(72,69){\circle*{2}}
\put(72,108){\circle*{2}}
\put(72,147){\circle*{2}}
\put(92,30){\circle*{2}}
\put(92,69){\circle*{2}}
\put(92,108){\circle*{2}}
\put(92,147){\circle*{2}}
\put(92,186){\circle*{2}}
\put(112,30){\circle*{2}}
\put(112,69){\circle*{2}}
\put(112,108){\circle*{2}}
\put(112,147){\circle*{2}}
\put(112,186){\circle*{2}}
\put(132,30){\circle*{2}}
\put(132,69){\circle*{2}}
\put(132,108){\circle*{2}}
\put(132,147){\circle*{2}}
\put(152,30){\circle*{2}}
\put(152,69){\circle*{2}}
\put(152,108){\circle*{2}}
\put(172,30){\circle*{2}}
\put(172,69){\circle*{2}}
\put(192,30){\circle*{2}}
\put(4,14){\circle*{2}}
\put(22,14){\circle*{2}}
\put(22,53){\circle*{2}}
\put(42,14){\circle*{2}}
\put(42,53){\circle*{2}}
\put(42,92){\circle*{2}}
\put(62,14){\circle*{2}}
\put(62,53){\circle*{2}}
\put(62,92){\circle*{2}}
\put(62,131){\circle*{2}}
\put(82,14){\circle*{2}}
\put(82,53){\circle*{2}}
\put(82,92){\circle*{2}}
\put(82,131){\circle*{2}}
\put(82,170){\circle*{2}}
\put(102,14){\circle*{2}}
\put(102,53){\circle*{2}}
\put(102,92){\circle*{2}}
\put(102,131){\circle*{2}}
\put(102,170){\circle*{2}}
\put(122,14){\circle*{2}}
\put(122,53){\circle*{2}}
\put(122,92){\circle*{2}}
\put(122,131){\circle*{2}}
\put(142,14){\circle*{2}}
\put(142,53){\circle*{2}}
\put(142,92){\circle*{2}}
\put(162,14){\circle*{2}}
\put(162,53){\circle*{2}}
\put(182,14){\circle*{2}}
\put(4,30){\circle*{2}}
\put(22,30){\circle*{2}}
\put(22,69){\circle*{2}}
\put(42,30){\circle*{2}}
\put(42,69){\circle*{2}}
\put(42,108){\circle*{2}}
\put(62,30){\circle*{2}}
\put(62,69){\circle*{2}}
\put(62,108){\circle*{2}}
\put(62,147){\circle*{2}}
\put(82,30){\circle*{2}}
\put(82,69){\circle*{2}}
\put(82,108){\circle*{2}}
\put(82,147){\circle*{2}}
\put(82,186){\circle*{2}}
\put(102,30){\circle*{2}}
\put(102,69){\circle*{2}}
\put(102,108){\circle*{2}}
\put(102,147){\circle*{2}}
\put(102,186){\circle*{2}}
\put(122,30){\circle*{2}}
\put(122,69){\circle*{2}}
\put(122,108){\circle*{2}}
\put(122,147){\circle*{2}}
\put(142,30){\circle*{2}}
\put(142,69){\circle*{2}}
\put(142,108){\circle*{2}}
\put(162,30){\circle*{2}}
\put(162,69){\circle*{2}}
\put(182,30){\circle*{2}}
\put(14,14){\line(0,1){0}}
\put(32,14){\line(0,1){0}}
\put(32,53){\line(0,1){0}}
\put(52,14){\line(0,1){0}}
\put(52,53){\line(0,1){0}}
\put(52,92){\line(0,1){0}}
\put(72,14){\line(0,1){0}}
\put(72,53){\line(0,1){0}}
\put(72,92){\line(0,1){0}}
\put(72,131){\line(0,1){0}}
\put(92,14){\line(0,1){0}}
\put(92,53){\line(0,1){0}}
\put(92,92){\line(0,1){0}}
\put(92,131){\line(0,1){0}}
\put(92,170){\line(0,1){0}}
\put(112,14){\line(0,1){0}}
\put(112,53){\line(0,1){0}}
\put(112,92){\line(0,1){0}}
\put(112,131){\line(0,1){0}}
\put(112,170){\line(0,1){0}}
\put(132,14){\line(0,1){0}}
\put(132,53){\line(0,1){0}}
\put(132,92){\line(0,1){0}}
\put(132,131){\line(0,1){0}}
\put(152,14){\line(0,1){0}}
\put(152,53){\line(0,1){0}}
\put(152,92){\line(0,1){0}}
\put(172,14){\line(0,1){0}}
\put(172,53){\line(0,1){0}}
\put(192,14){\line(0,1){0}}
\put(14,30){\line(0,1){0}}
\put(32,30){\line(0,1){0}}
\put(32,69){\line(0,1){0}}
\put(52,30){\line(0,1){0}}
\put(52,69){\line(0,1){0}}
\put(52,108){\line(0,1){0}}
\put(72,30){\line(0,1){0}}
\put(72,69){\line(0,1){0}}
\put(72,108){\line(0,1){0}}
\put(72,147){\line(0,1){0}}
\put(92,30){\line(0,1){0}}
\put(92,69){\line(0,1){0}}
\put(92,108){\line(0,1){0}}
\put(92,147){\line(0,1){0}}
\put(92,186){\line(0,1){0}}
\put(112,30){\line(0,1){0}}
\put(112,69){\line(0,1){0}}
\put(112,108){\line(0,1){0}}
\put(112,147){\line(0,1){0}}
\put(112,186){\line(0,1){0}}
\put(132,30){\line(0,1){0}}
\put(132,69){\line(0,1){0}}
\put(132,108){\line(0,1){0}}
\put(132,147){\line(0,1){0}}
\put(152,30){\line(0,1){0}}
\put(152,69){\line(0,1){0}}
\put(152,108){\line(0,1){0}}
\put(172,30){\line(0,1){0}}
\put(172,69){\line(0,1){0}}
\put(192,30){\line(0,1){0}}
\put(2,2){\framebox(14,14)[]{}}
\put(20,2){\framebox(14,14)[]{}}
\put(20,41){\framebox(14,14)[]{}}
\put(40,2){\framebox(14,14)[]{}}
\put(40,41){\framebox(14,14)[]{}}
\put(40,80){\framebox(14,14)[]{}}
\put(60,2){\framebox(14,14)[]{}}
\put(60,41){\framebox(14,14)[]{}}
\put(60,80){\framebox(14,14)[]{}}
\put(60,119){\framebox(14,14)[]{}}
\put(80,2){\framebox(14,14)[]{}}
\put(80,41){\framebox(14,14)[]{}}
\put(80,80){\framebox(14,14)[]{}}
\put(80,119){\framebox(14,14)[]{}}
\put(80,158){\framebox(14,14)[]{}}
\put(100,2){\framebox(14,14)[]{}}
\put(100,41){\framebox(14,14)[]{}}
\put(100,80){\framebox(14,14)[]{}}
\put(100,119){\framebox(14,14)[]{}}
\put(100,158){\framebox(14,14)[]{}}
\put(120,2){\framebox(14,14)[]{}}
\put(120,41){\framebox(14,14)[]{}}
\put(120,80){\framebox(14,14)[]{}}
\put(120,119){\framebox(14,14)[]{}}
\put(140,2){\framebox(14,14)[]{}}
\put(140,41){\framebox(14,14)[]{}}
\put(140,80){\framebox(14,14)[]{}}
\put(160,2){\framebox(14,14)[]{}}
\put(160,41){\framebox(14,14)[]{}}
\put(180,2){\framebox(14,14)[]{}}
\put(2,18){\framebox(14,14)[]{}}
\put(20,18){\framebox(14,14)[]{}}
\put(20,57){\framebox(14,14)[]{}}
\put(40,18){\framebox(14,14)[]{}}
\put(40,57){\framebox(14,14)[]{}}
\put(40,96){\framebox(14,14)[]{}}
\put(60,18){\framebox(14,14)[]{}}
\put(60,57){\framebox(14,14)[]{}}
\put(60,96){\framebox(14,14)[]{}}
\put(60,135){\framebox(14,14)[]{}}
\put(80,18){\framebox(14,14)[]{}}
\put(80,57){\framebox(14,14)[]{}}
\put(80,96){\framebox(14,14)[]{}}
\put(80,135){\framebox(14,14)[]{}}
\put(80,174){\framebox(14,14)[]{}}
\put(100,18){\framebox(14,14)[]{}}
\put(100,57){\framebox(14,14)[]{}}
\put(100,96){\framebox(14,14)[]{}}
\put(100,135){\framebox(14,14)[]{}}
\put(100,174){\framebox(14,14)[]{}}
\put(120,18){\framebox(14,14)[]{}}
\put(120,57){\framebox(14,14)[]{}}
\put(120,96){\framebox(14,14)[]{}}
\put(120,135){\framebox(14,14)[]{}}
\put(140,18){\framebox(14,14)[]{}}
\put(140,57){\framebox(14,14)[]{}}
\put(140,96){\framebox(14,14)[]{}}
\put(160,18){\framebox(14,14)[]{}}
\put(160,57){\framebox(14,14)[]{}}
\put(180,18){\framebox(14,14)[]{}}
\put(4,30){\line(1,0){10}}
\put(82,30){\line(1,0){10}}
\put(82,69){\line(1,0){10}}
\put(82,108){\line(1,0){10}}
\put(82,147){\line(1,0){10}}
\put(82,186){\line(1,0){10}}
\put(122,30){\line(1,0){10}}
\put(122,69){\line(1,0){10}}
\put(122,108){\line(1,0){10}}
\put(122,147){\line(1,0){10}}
\put(142,30){\line(1,0){10}}
\put(142,69){\line(1,0){10}}
\put(142,108){\line(1,0){10}}
\put(162,30){\line(1,0){10}}
\put(162,69){\line(1,0){10}}
\put(182,30){\line(1,0){10}}

{\footnotesize

\put(9,34.5){\makebox(0,0)[cc]{P (\ref{observation1})}}%10
\put(27,34.5){\makebox(0,0)[cc]{NPC (\ref{theorem11})}}%11
\put(47,34.5){\makebox(0,0)[cc]{NPC (\ref{theorem12})}}%12
\put(67,34.5){\makebox(0,0)[cc]{}}%13
\put(87,34.5){\makebox(0,0)[cc]{}}%14
\put(107,34.5){\makebox(0,0)[cc]{}}%15
\put(127,34.5){\makebox(0,0)[cc]{NPC (\ref{theorem16})}}%16
\put(147,34.5){\makebox(0,0)[cc]{}}%17
\put(167,34.5){\makebox(0,0)[cc]{}}%18
\put(187,34.5){\makebox(0,0)[cc]{P (\ref{theorem19})}}%19

\put(27,73.5){\makebox(0,0)[cc]{P (\ref{observation1})}}%20
\put(47,73.5){\makebox(0,0)[cc]{P (\ref{observation1})}}%21
\put(67,73.5){\makebox(0,0)[cc]{}}%22
\put(87,73.5){\makebox(0,0)[cc]{P (\ref{theorem23})}}%23
\put(107,73.5){\makebox(0,0)[cc]{}}%24
\put(127,73.5){\makebox(0,0)[cc]{NPC (\ref{theorem25})}}%25
\put(147,73.5){\makebox(0,0)[cc]{NPC (\ref{theorem26})}}%26
\put(167,73.5){\makebox(0,0)[cc]{NPC(\ref{theorem27})}}%27

\put(47,112.5){\makebox(0,0)[cc]{P (\ref{theorem30})}}%30
\put(67,112.5){\makebox(0,0)[cc]{}}%31
\put(87,112.5){\makebox(0,0)[cc]{P (\ref{theorem32})}}%32
\put(107,112.5){\makebox(0,0)[cc]{P (\ref{theorem33})}}%33
\put(127,112.5){\makebox(0,0)[cc]{NPC (\ref{theorem34})}}%34
\put(147,112.5){\makebox(0,0)[cc]{P (\ref{theorem35})}}%35

\put(67,151.5){\makebox(0,0)[cc]{P (\ref{theorem40})}}%40
\put(87,151.5){\makebox(0,0)[cc]{P (\ref{theorem41+50})}}%41
\put(107,151.5){\makebox(0,0)[cc]{P (\ref{theorem42})}}%42
\put(127,151.5){\makebox(0,0)[cc]{P (\ref{theorem43})}}%43

\put(87,190.5){\makebox(0,0)[cc]{P (\ref{theorem41+50})}}%50
\put(107,190.5){\makebox(0,0)[cc]{P (\ref{theorem52})}}%52

}

\put(0,0){\framebox(18,37)[cc]{}}
\put(18,0){\framebox(18,37)[cc]{}}
\put(18,39){\framebox(18,37)[cc]{}}
\put(38,0){\framebox(18,37)[cc]{}}
\put(38,39){\framebox(18,37)[cc]{}}
\put(38,78){\framebox(18,37)[cc]{}}
\put(58,0){\framebox(18,37)[cc]{}}
\put(58,39){\framebox(18,37)[cc]{}}
\put(58,78){\framebox(18,37)[cc]{}}
\put(58,117){\framebox(18,37)[cc]{}}
\put(78,0){\framebox(18,37)[cc]{}}
\put(78,39){\framebox(18,37)[cc]{}}
\put(78,78){\framebox(18,37)[cc]{}}
\put(78,117){\framebox(18,37)[cc]{}}
\put(78,156){\framebox(18,37)[cc]{}}
\put(98,0){\framebox(18,37)[cc]{}}
\put(98,39){\framebox(18,37)[cc]{}}
\put(98,78){\framebox(18,37)[cc]{}}
\put(98,117){\framebox(18,37)[cc]{}}
\put(98,156){\framebox(18,37)[cc]{}}
\put(118,0){\framebox(18,37)[cc]{}}
\put(118,39){\framebox(18,37)[cc]{}}
\put(118,78){\framebox(18,37)[cc]{}}
\put(118,117){\framebox(18,37)[cc]{}}
\put(138,0){\framebox(18,37)[cc]{}}
\put(138,39){\framebox(18,37)[cc]{}}
\put(138,78){\framebox(18,37)[cc]{}}
\put(158,0){\framebox(18,37)[cc]{}}
\put(158,39){\framebox(18,37)[cc]{}}
\put(178,0){\framebox(18,37)[cc]{}}
\put(22,53){\line(1,0){10}}
\put(42,53){\line(1,0){10}}
\put(62,53){\line(1,0){10}}
\put(82,53){\line(1,0){10}}
\put(82,170){\line(1,0){10}}
\put(102,53){\line(1,0){10}}
\put(102,170){\line(1,0){10}}
\put(122,53){\line(1,0){10}}
\put(142,53){\line(1,0){10}}
\put(162,53){\line(1,0){10}}
\put(22,69){\line(0,-1){10}}
\put(22,59){\line(1,1){10}}
\put(22,30){\line(0,-1){10}}
\put(22,20){\line(1,1){10}}
\put(42,30){\line(0,-1){11}}
\put(42,108){\line(0,-1){11}}
\put(52,30){\line(0,-1){10}}
\put(52,69){\line(0,-1){10}}
\put(52,108){\line(0,-1){10}}
\put(72,30){\line(-1,0){11}}
\put(72,69){\line(-1,0){11}}
\put(72,108){\line(-1,0){11}}
\put(72,147){\line(-1,0){11}}
\put(62,30){\line(0,-1){10}}
\put(62,69){\line(0,-1){10}}
\put(62,108){\line(0,-1){10}}
\put(62,147){\line(0,-1){10}}
\put(62,20){\line(1,1){10}}
\put(62,59){\line(1,1){10}}
\put(62,98){\line(1,1){10}}
\put(62,137){\line(1,1){10}}
\put(42,92){\line(0,-1){10}}
\put(42,82){\line(1,1){10}}
\put(42,69){\line(0,-1){10}}
\put(62,92){\line(0,-1){10}}
\put(82,92){\line(0,-1){10}}
\put(102,92){\line(0,-1){10}}
\put(122,92){\line(0,-1){10}}
\put(142,92){\line(0,-1){10}}
\put(72,92){\line(0,-1){10}}
\put(92,92){\line(0,-1){10}}
\put(112,92){\line(0,-1){10}}
\put(132,92){\line(0,-1){10}}
\put(152,92){\line(0,-1){10}}
\put(62,131){\line(0,-1){10}}
\put(82,131){\line(0,-1){10}}
\put(82,170){\line(0,-1){10}}
\put(102,131){\line(0,-1){10}}
\put(102,170){\line(0,-1){10}}
\put(122,131){\line(0,-1){10}}
\put(62,121){\line(1,1){10}}
\put(82,121){\line(1,1){10}}
\put(82,160){\line(1,1){10}}
\put(102,121){\line(1,1){10}}
\put(102,160){\line(1,1){10}}
\put(122,121){\line(1,1){10}}
\put(82,20){\line(0,1){10}}
\put(82,59){\line(0,1){10}}
\put(82,98){\line(0,1){10}}
\put(82,137){\line(0,1){10}}
\put(82,176){\line(0,1){10}}
\put(92,20){\line(0,1){10}}
\put(92,59){\line(0,1){10}}
\put(92,98){\line(0,1){10}}
\put(92,137){\line(0,1){10}}
\put(92,176){\line(0,1){10}}
\put(112,20){\line(-1,0){10}}
\put(112,59){\line(-1,0){10}}
\put(112,98){\line(-1,0){10}}
\put(112,137){\line(-1,0){10}}
\put(112,176){\line(-1,0){10}}
\put(102,20){\line(1,1){10}}
\put(102,59){\line(1,1){10}}
\put(102,98){\line(1,1){10}}
\put(102,137){\line(1,1){10}}
\put(102,176){\line(1,1){10}}
\put(102,30){\line(0,-1){10}}
\put(102,69){\line(0,-1){10}}
\put(102,108){\line(0,-1){10}}
\put(102,147){\line(0,-1){10}}
\put(102,186){\line(0,-1){10}}
\put(122,30){\line(0,-1){10}}
\put(122,69){\line(0,-1){10}}
\put(122,108){\line(0,-1){10}}
\put(122,147){\line(0,-1){10}}
\put(122,20){\line(1,0){10}}
\put(122,59){\line(1,0){10}}
\put(122,98){\line(1,0){10}}
\put(122,137){\line(1,0){10}}
\put(132,30){\line(-1,-1){10}}
\put(132,69){\line(-1,-1){10}}
\put(132,108){\line(-1,-1){10}}
\put(132,147){\line(-1,-1){10}}
\put(142,30){\line(0,-1){10}}
\put(142,69){\line(0,-1){10}}
\put(142,108){\line(0,-1){10}}
\put(142,20){\line(1,0){10}}
\put(142,59){\line(1,0){10}}
\put(142,98){\line(1,0){10}}
\put(152,20){\line(0,1){10}}
\put(152,59){\line(0,1){10}}
\put(152,98){\line(0,1){10}}
\put(162,69){\line(0,-1){10}}
\put(162,59){\line(1,0){10}}
\put(172,59){\line(0,1){10}}
\put(172,69){\line(-1,-1){10}}
\put(162,30){\line(0,-1){10}}
\put(182,30){\line(0,-1){10}}
\put(162,20){\line(1,0){10}}
\put(182,20){\line(1,0){10}}
\put(172,20){\line(0,1){10}}
\put(192,20){\line(0,1){10}}
\put(172,30){\line(-1,-1){10}}
\put(192,30){\line(-1,-1){10}}
\put(182,30){\line(1,-1){10}}
\end{picture}
\end{center}
\caption{All $30$ pairs of non-isomorphic graphs of order four up to complementation,
together with the status of the corresponding sandwich decision problem,
where ``P'' means ``{\it polynomial time solvable}'',
``NPC'' means ``{\it NP-complete}'', and
the number in the bracket is the reference number of the corresponding result within this paper.}\label{fig1}
\end{figure}
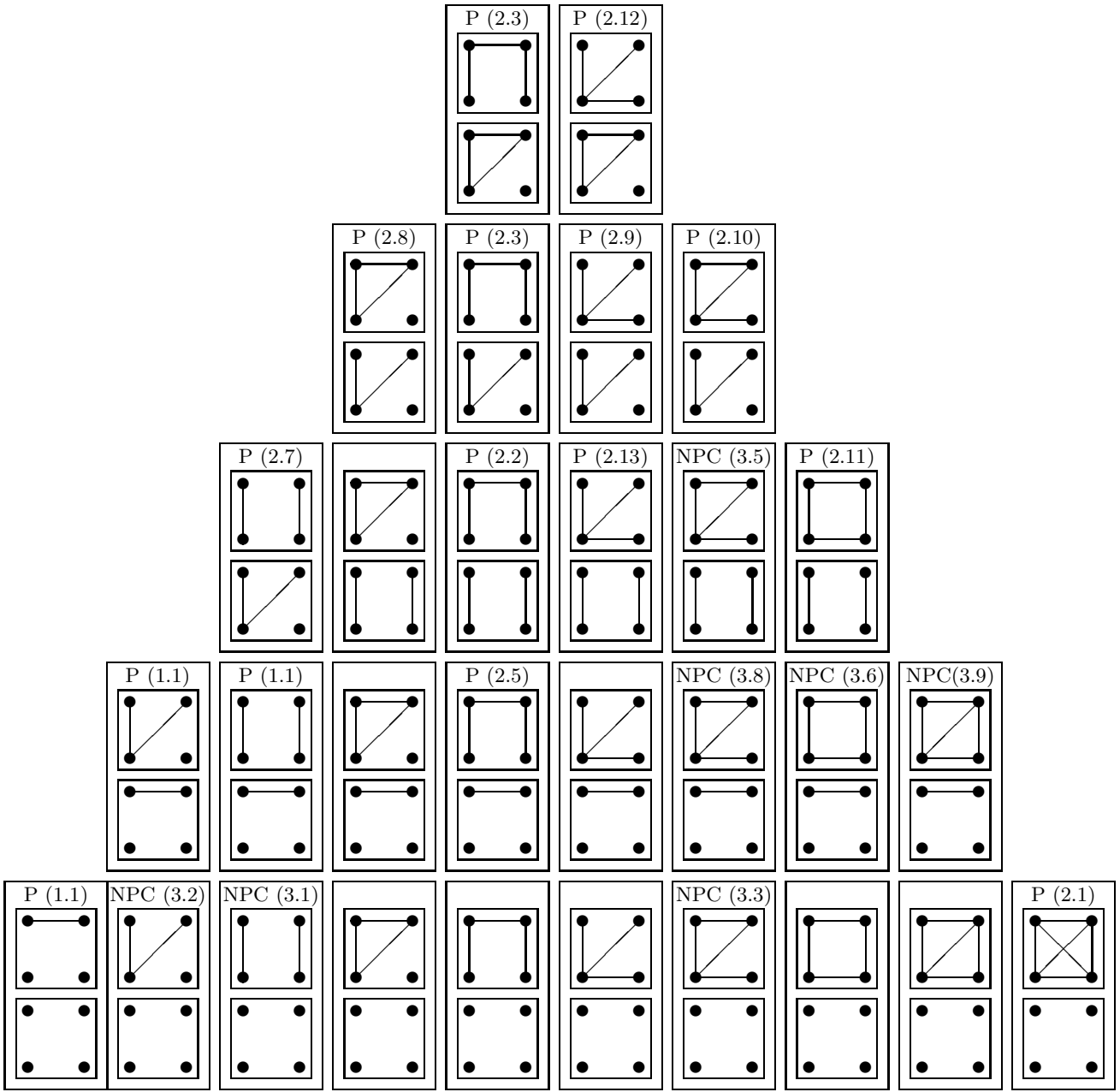

\noindent In Sections 2 and 3, we collect our positive and negative results, respectively.
In a final section, we conclude with some comments on the open cases. 

\section{Some Tractable Cases}

We present our positive results in an order of roughly increasing difficulty.

Observation \ref{observation1}(i) and (iv) imply that the ${\cal F}$-{\sc free Sandwich Problem} can be solved in polynomial time 
if ${\cal F}$ or $\overline{\cal F}$ is one of the sets
$$\{ P_3\},
\left\{ {\rm diamond},K_4\right\},
\left\{ {\rm diamond},C_4\right\},\mbox{ or }
\left\{ {\rm diamond},{\rm paw}\right\}.$$
In order to understand the complexity of the ${\cal F}$-{\sc free Sandwich Problem},
if ${\cal F}$ is as in Observation \ref{observation1}(ii), then it suffices to consider instances $(G_1,G_2)$ such that $G_1$ is connected,
and,
if ${\cal F}$ is as in Observation \ref{observation1}(iii), then it suffices to consider instances $(G_1,G_2)$ such that $G_2$ has no universal vertex;
otherwise, in both cases some simple algorithmic reduction applies.

For positive integers $s$ and $t$, 
let $R(s,t)$ be the {\it Ramsey number},
in particular, every graph of order $R(s,t)$ contains an induced $K_s$ or $\overline{K_r}$.

\begin{theorem}\label{theorem19}
The $\left\{ K_4,\overline{K_4}\right\}$-{\sc free Sandwich Problem} can be solved in polynomial time.
\end{theorem}
{\it Proof:} No instance $(G_1,G_2)$ of the $\left\{ K_4,\overline{K_4}\right\}$-{\sc free Sandwich Problem}
for which $G_1$ has order at least $R(4,4)$ has a solution.
Instances $(G_1,G_2)$ for which $G_1$ has order less than $R(4,4)$ can be solved in constant time. $\Box$

\medskip

\noindent It is well-known \cite{brlesp} that, for every $P_4$-free graph $G$ or order at least $2$, either $G$ or $\overline{G}$ is disconnected.

\begin{theorem}\label{theorem32}
The $\{ P_4,C_4\}$-{\sc free Sandwich Problem} can be solved in polynomial time.
\end{theorem}
{\it Proof:} Let $(G_1,G_2)$ be an instance of the ${\cal F}$-{\sc free Sandwich Problem} for ${\cal F}=\{ P_4,C_4\}$.
By Observation \ref{observation1}(ii) and (iii), we may assume that $G_1$ is connected, and that $G_2$ has no universal vertex.
Suppose that ${\cal SW}_{\cal F}(G_1,G_2)$ contains some graph $G$.
Let $u$ is a vertex of maximum degree in $G$. 
Since $u$ is not universal, there is an induced path $uvw$ of order $3$ in $G$.
Since $w$ is a neighbor of $v$ but not of $u$, and $u$ has at least as many neighbors as $v$,
there is a vertex $x$ that is a neighbor of $u$ but not of $v$.
Nevertheless, the subgraph $G[\{ u,v,w,x\}]$ of $G$ induced by $\{ u,v,w,x\}$ is either $P_4$ or $C_4$, which is a contradiction.
Hence, either one of the two algorithmic reductions corresponding to Observation \ref{observation1}(ii) and (iii) applies,
or ${\cal SW}_{\cal F}(G_1,G_2)$ is necessarily empty. $\Box$

\medskip

\noindent The following result concerns the two cases 
${\cal F}=\left\{ P_4,\overline{\rm paw}\right\}$ and ${\cal F}=\left\{ P_4,\overline{\rm claw}\right\}$.

\begin{theorem}\label{theorem41+50}
The $\{ P_4,K_1\cup F_2\}$-{\sc free Sandwich Problem} can be solved in polynomial time for $F_2\in \{ K_3,P_3\}$.
\end{theorem}
{\it Proof:} Let $(G_1,G_2)$ be an instance of the ${\cal F}$-{\sc free Sandwich Problem} for ${\cal F}=\{ P_4,K_1\cup F_2\}$,
where $n(G_1)\geq 2$.
If ${\cal SW}_{\cal F}(G_1,G_2)$ contains some disconnected graph $G$,
then $G$ is $\{ P_4,F_2\}$-free.
Since, by Observation \ref{observation1}(iv), 
the $\{ P_4,F_2\}$-{\sc free Sandwich Problem} can be solved in polynomial time,
this possibility can be checked in polynomial time.
Hence, we may assume that ${\cal SW}_{\cal F}(G_1,G_2)$ contains no disconnected graph.
If $\overline{G_2}$ is connected, then, since, for every graph $G$ in ${\cal SW}_{\cal F}(G_1,G_2)$,
the graph $\overline{G}$ is $P_4$-free and contains $\overline{G_2}$,
${\cal SW}_{\cal F}(G_1,G_2)$ is empty.
Hence, we may assume that $\overline{G_2}$ is disconnected.
Note that, if $H$ is the join of two graphs $H_1$ and $H_2$, then every induced $P_4$ or $K_1\cup F_2$ in $H$
is completely contained either in $H_1$ or in $H_2$.
Hence, if $K$ is the vertex set of some component of $\overline{G_2}$,
then 
${\cal SW}_{\cal F}(G_1,G_2)$ is non-empty
if and only if 
${\cal SW}_{\cal F}(G_1[K],G_2[K])$ and
${\cal SW}_{\cal F}(G_1-K,G_2-K)$ are both non-empty,
that is, in polynomial time, one can reduce the instance $(G_1,G_2)$
to two smaller instances $(G'_1,G'_2)$ and $(G''_1,G''_2)$ such that $n(G_1)=n(G_1')+n(G_1'')$,
which implies the desired statement. $\Box$

\medskip

\noindent Deciding the existence of a complete bipartite sandwich can easily be reduced to $2${\sc Sat} \cite{gokash,tedafi}.
We give a different argument leading to a simpler algorithm.

\begin{lemma}\label{lemma_complete_bipartite}
If $\Pi$ is the set of all complete bipartite graphs, 
then the $\Pi$-{\sc Sandwich Problem} can be solved in polynomial time.
\end{lemma}
{\it Proof:} Let $(G_1,G_2)$ be an instance of the $\Pi$-{\sc Sandwich Problem}.
Clearly, we may assume that all components $K_1,\ldots,K_p$ of $G_1$ are bipartite.
Let $K_i$ have the partite sets $A_i$ and $B_i$ for $i\in [p]$.
Initialize a set ${\cal P}$ as $\{ \{ A_1,B_1\},\ldots,\{ A_p,B_p\}\}$, 
and, iteratively and as long as possible, 
whenever ${\cal P}$ contains two distinct sets $\{ X,Y\}$ and $\{ X',Y'\}$ 
such that, in $G_2$, some vertex in $X$ is non-adjacent to some vertex in $X'$,
then replace $\{ X,Y\}$ and $\{ X',Y'\}$ within ${\cal P}$ by $\{ X\cup X',Y\cup Y'\}$;
breaking ties arbitrarily.
Note that $X\cup X'$ is the union of partite sets of components of $G_1$
that necessarily belong to the same partite set of any solution.
When ${\cal P}$ no longer changes, then, 
for every two distinct sets $\{ X,Y\}$ and $\{ X',Y'\}$ in ${\cal P}$,
the graph $G_2$ contains all edges between $X\cup Y$ and $X'\cup Y'$.
Therefore, if the final ${\cal P}$ contains the sets $\{ X_1,Y_1\},\ldots,\{ X_q,Y_q\}$,
then there is a complete bipartite graph $G$ with $G_1\subseteq G\subseteq G_2$
if and only if $G_2$ contains all edges between $X_i$ and $Y_i$ for every $i\in [q]$.
Furthermore, such a graph $G$ can easily be determined.
$\Box$

\begin{theorem}\label{theorem23}
The $\left\{ P_4,\overline{\rm diamond}\right\}$-{\sc free Sandwich Problem} can be solved in polynomial time.
\end{theorem}
{\it Proof:} Let $(G_1,G_2)$ be an instance of the ${\cal F}$-{\sc free Sandwich Problem} for ${\cal F}=\left\{ P_4,\overline{\rm diamond}\right\}$,
where $m(G_1)>0$.

Suppose that ${\cal SW}_{\cal F}(G_1,G_2)$ contains some disconnected graph $G$.
Since $G$ is $P_4$-free, 
$\overline{G}$ is a connected graph in ${\cal SW}_{\overline{\cal F}}\left(\overline{G_2},\overline{G_1}\right)$.
Let $\overline{G}$ be the join of the two non-empty graphs $\overline{G}_L$ and $\overline{G}_R$.
Since $\overline{G}$ is diamond-free, 
the two graphs $\overline{G}_L$ and $\overline{G}_R$ are $P_3$-free,
that is, they are the unions of $k_L$ and $k_R$ complete graphs, respectively.
Since $G$ has at least one edge, we may assume, by symmetry, that $k_R\geq 2$.
Since $\overline{G}$ is diamond-free, this implies that all vertices of $\overline{G}_L$ are isolated.
If $k_L=1$, then $\overline{G}_L$ consists of a universal vertex $u_L$ of $\overline{G}$. 
Since the $\{ P_3\}$-{\sc free Sandwich Problem} can be solved in polynomial time,
considering all $n(G_1)$ choices for $u_L$, 
one can check in polynomial time whether ${\cal SW}_{\cal F}(G_1,G_2)$ contains such a graph $G$.
Hence, we may assume that $k_L\geq 2$.
Since $\overline{G}$ is diamond-free, this implies that all vertices of $\overline{G}_R$ are isolated,
that is, $\overline{G}$ is a complete bipartite graph.
By Lemma \ref{lemma_complete_bipartite},
one can check in polynomial time 
whether ${\cal SW}_{\overline{\cal F}}\left(\overline{G_2},\overline{G_1}\right)$ 
contains a complete bipartite graph.
Altogether, it follows that one can check in polynomial time 
whether ${\cal SW}_{\cal F}(G_1,G_2)$ contains some disconnected graph.
Hence, we may assume that ${\cal SW}_{\cal F}(G_1,G_2)$ contains no disconnected graph.

If $\overline{G_2}$ is connected, then, similarly as in the proof of Theorem \ref{theorem41+50},
${\cal SW}_{\cal F}(G_1,G_2)$ is empty.
Hence, we may assume that $\overline{G_2}$ is disconnected.
If $K$ is the vertex set of some component of $\overline{G_2}$,
then, similarly as in the proof of Theorem \ref{theorem41+50}, 
${\cal SW}_{\cal F}(G_1,G_2)$ is non-empty
if and only if 
${\cal SW}_{\cal F}(G_1[K],G_2[K])$ and
${\cal SW}_{\cal F}(G_1-K,G_2-K)$ are both non-empty,
which implies the desired statement. $\Box$

\medskip

\noindent Our next few results involve the paw, and the following result of Olariu is quite useful.

\begin{lemma}[Olariu \cite{ol}]\label{lemma_olariu}
A connected graph is ${\rm paw}$-free if and only if it is triangle-free or $\overline{P_3}$-free.
\end{lemma}

\begin{theorem}\label{theorem30}
The $\left\{ {\rm paw},C_4\right\}$-{\sc free Sandwich Problem} can be solved in polynomial time.
\end{theorem}
{\it Proof:} Let $(G_1,G_2)$ be an instance of the ${\cal F}$-{\sc free Sandwich Problem} 
for ${\cal F}=\left\{ {\rm paw},C_4\right\}$.
By Observation \ref{observation1}(ii), we may assume that $G_1$ is connected.

Suppose that ${\cal SW}_{\cal F}(G_1,G_2)$ contains some graph $G$.
By Lemma \ref{lemma_olariu}, $G$ is triangle-free or $\overline{P_3}$-free.
Since $G$ is $\left\{ K_3,C_4\right\}$-free if and only if $G=G_1$, and $G_1$ is $\left\{ K_3,C_4\right\}$-free,
we may assume that $G$ is $\left\{ \overline{P_3},C_4\right\}$-free.
This implies that $G$ is a complete multipartite graph with at most one partite set of order more than $1$,
that is, $\overline{G}$ has at most one edge.
This implies that $\overline{G_2}$ has at most one edge, 
and, hence, that $G_2$ is $\left\{ \overline{P_3},C_4\right\}$-free.
Altogether, if ${\cal SW}_{\cal F}(G_1,G_2)$ is non-empty,
then $G_1$ or $G_2$ belongs to this set.
$\Box$

\begin{theorem}\label{theorem40}
The $\left\{ {\rm paw},{\rm claw}\right\}$-{\sc free Sandwich Problem} can be solved in polynomial time.
\end{theorem}
{\it Proof:} Let $(G_1,G_2)$ be an instance of the ${\cal F}$-{\sc free Sandwich Problem} 
for ${\cal F}=\left\{ {\rm paw},{\rm claw}\right\}$.
By Observation \ref{observation1}(ii), we may assume that $G_1$ is connected.

Suppose that ${\cal SW}_{\cal F}(G_1,G_2)$ contains some graph $G$.
By Lemma \ref{lemma_olariu}, $G$ is triangle-free or $\overline{P_3}$-free.
Since $G$ is $\left\{ K_3,{\rm claw}\right\}$-free if and only if $G=G_1$, and $G_1$ is $\left\{ K_3,{\rm claw}\right\}$-free,
we may assume that $G$ is $\left\{ \overline{P_3},{\rm claw}\right\}$-free.
This implies that $G$ is a complete multipartite graph such that each partite set contains at most two vertices,
that is, $\overline{G}$ has maximum degree at most $1$.
This implies that $\overline{G_2}$ has maximum degree at most $1$, 
and, hence, that $G_2$ is $\left\{ \overline{P_3},{\rm claw}\right\}$-free.
Similarly as in the proof of Theorem \ref{theorem30}, if ${\cal SW}_{\cal F}(G_1,G_2)$ is non-empty,
then $G_1$ or $G_2$ belong to this set.
$\Box$

\begin{theorem}\label{theorem42}
The $\left\{ {\rm paw},\overline{{\rm claw}}\right\}$-{\sc free Sandwich Problem} can be solved in polynomial time.
\end{theorem}
{\it Proof:} Let $(G_1,G_2)$ be an instance of the ${\cal F}$-{\sc free Sandwich Problem} 
for ${\cal F}=\left\{ {\rm paw},\overline{{\rm claw}}\right\}$.

Suppose that ${\cal SW}_{\cal F}(G_1,G_2)$ contains some graph $G$.
If $G$ is triangle-free, then $G_1$ is also triangle-free, and, hence, lies in ${\cal SW}_{\cal F}(G_1,G_2)$.
Hence, we may assume that $G$ is not triangle-free.
Since $G$ is $\overline{{\rm claw}}$-free, $G$ is connected.
By Lemma \ref{lemma_olariu}, $G$ is $\overline{P_3}$-free.
Since a graph is $\left\{ \overline{P_3},\overline{\rm claw}\right\}$-free 
if and only if it is $\overline{P_3}$-free,
and the $\left\{ \overline{P_3}\right\}$-{\sc free Sandwich Problem} can be solved in polynomial time,
one can check in polynomial time whether ${\cal SW}_{\cal F}(G_1,G_2)$ contains such a graph.
$\Box$

\begin{theorem}\label{theorem43}
The $\left\{ {\rm paw},\overline{\rm paw}\right\}$-{\sc free Sandwich Problem} can be solved in polynomial time.
\end{theorem}
{\it Proof:} Let $(G_1,G_2)$ be an instance of the ${\cal F}$-{\sc free Sandwich Problem} 
for ${\cal F}=\left\{ {\rm paw},\overline{\rm paw}\right\}$.

Suppose that ${\cal SW}_{\cal F}(G_1,G_2)$ contains some graph $G$.
Since $\overline{\cal F}={\cal F}$, we may assume, by Observation \ref{observation1}, that $G$ is connected.
By Lemma \ref{lemma_olariu}, this implies that $G$ is triangle-free or $\overline{P_3}$-free.
Since a graph is $\left\{ \overline{\rm paw},\overline{P_3}\right\}$-free if and only if it is $\overline{P_3}$-free,
and the $\left\{ \overline{P_3}\right\}$-{\sc free Sandwich Problem} can be solved in polynomial time,
one can check in polynomial time 
whether ${\cal SW}_{\cal F}(G_1,G_2)$ contains a $\left\{ \overline{\rm paw},\overline{P_3}\right\}$-free graph.
Hence, we may assume that $G$ is $\left\{ \overline{\rm paw},K_3\right\}$-free.

By Lemma \ref{lemma_complete_bipartite}, we may assume that $G$ is not complete bipartite.
If the maximum degree of $G$ is at most $2$, then $G$ has at most $5$ vertices.
Hence, we may assume that $G$ has maximum degree at least $3$.
Let $u$ be a vertex of maximum degree.
Let $B$ be the neighborhood of $u$.
Let $A$ be the set of vertices whose neighborhood is $B$.
Since $G$ is triangle-free, $G[A\cup B]$ is a complete bipartite graph with partite sets $A$ and $B$.
Since $G$ is connected but not complete bipartite,
some vertex in $B$ has a neighbor $w$ outside of $A$.
By the definition of $A$, $w$ has a non-neighbor $v$ in $B$.
Since $G$ is $\overline{\rm paw}$-free, $v$ is the only non-neighbor of $w$ in $B$.
Now, $v$, $w$, and two further vertices from $B$ induce a $\overline{\rm paw}$,
which is a contradiction. 
$\Box$

\medskip

\noindent Our next result relies on Maffray and Preissmann's \cite{mapr}
characterization of pseudo-split graphs,
and Golumbic, Kaplan, and Shamir's \cite{gokash} algorithm 
for the split sandwich problem.

\begin{theorem}\label{theorem35}
The $\left\{ C_4,\overline{C_4}\right\}$-{\sc free Sandwich Problem} can be solved in polynomial time.
\end{theorem}
{\it Proof:} Let $(G_1,G_2)$ be an instance of the ${\cal F}$-{\sc free Sandwich Problem} 
for ${\cal F}=\left\{ C_4,\overline{C_4}\right\}$.
Suppose that ${\cal SW}_{\cal F}(G_1,G_2)$ contains some graph $G$.
By a result of Maffray and Preissmann \cite{mapr},
there is a set $C$ of at most five vertices such that $G-C$ is a split graph.
Considering the $O\left(n(G_1)^5\right)$ choices for $C$,
and applying the polynomial time algorithm of Golumbic, Kaplan, and Shamir \cite{gokash}
to $G-C$, one can decide in polynomial time whether ${\cal SW}_{\cal F}(G_1,G_2)$ is non-empty. $\Box$

\medskip

\noindent The next proof uses
a result of 
Brandst\"{a}dt and Mahfud \cite{brma} 
concerning prime $\left\{ {\rm claw},\overline{{\rm claw}}\right\}$-free graphs.

\begin{theorem}\label{theorem52}
The $\left\{ {\rm claw},\overline{{\rm claw}}\right\}$-{\sc free Sandwich Problem} can be solved in polynomial time.
\end{theorem}
{\it Proof:} Let $(G_1,G_2)$ be an instance of the ${\cal F}$-{\sc free Sandwich Problem} 
for ${\cal F}=\left\{ {\rm claw},\overline{{\rm claw}}\right\}$, with $n(G_1)\geq 10$.
Suppose that ${\cal SW}_{\cal F}(G_1,G_2)$ contains some graph $G$.
We will show that either $G$ or $\overline{G}$ has maximum degree at most $2$,
that is, it is the union of paths and cycles.
If $G$ has maximum degree at most $2$, then $G_1$ belong to ${\cal SW}_{\cal F}(G_1,G_2)$, and,
if $\overline{G}$ has maximum degree at most $2$, then $G_2$ belong to ${\cal SW}_{\cal F}(G_1,G_2)$,
which clearly implies the desired statement.

If $G$ is disconnected, then $G$ is $\left\{ {\rm claw},K_3\right\}$-free, 
which clearly implies that $G$ has maximum degree at most $2$.
Hence, by symmetry, we may assume that $G$ and $\overline{G}$ are both connected.
If $G$ is prime, then Brandst\"{a}dt and Mahfud \cite{brma} 
showed that $G$ or $\overline{G}$ has maximum degree at most $2$.
Hence, we may assume that $G$ contains a homogeneous set $U$ of vertices,
that is, $2\leq |U|\leq n(G)-1$, and $V(G)=U\cup A\cup B$,
where $A$ is the set of vertices in $V(G)\setminus U$ that are adjacent to every vertex in $U$, and,
$B$ is the set of vertices in $V(G)\setminus U$ that are adjacent to no vertex in $U$.
Since $G$ and $\overline{G}$ are both connected, 
there are vertices $a$ and $a'$ in $A$, and, $b$ and $b'$ in $B$ such that 
$a$ and $b$ are adjacent,
and, 
$a'$ and $b'$ are non-adjacent.
Let $u$ and $u'$ be two vertices in $U$.
If $u$ and $u'$ are adjacent, then $G[\{ u,u',a',b'\}]$ is a $\overline{{\rm claw}}$, and,
if $u$ and $u'$ are not adjacent, then $G[\{ u,u',a,b\}]$ is a claw,
which completes the proof. $\Box$

\begin{theorem}\label{theorem33}
The $\left\{ {\rm claw},\overline{C_4}\right\}$-{\sc free Sandwich Problem} can be solved in polynomial time.
\end{theorem}
{\it Proof:} Let $(G_1,G_2)$ be an instance of the ${\cal F}$-{\sc free Sandwich Problem} for ${\cal F}=\left\{ \overline{\rm claw},C_4\right\}$,
which, by Observation \ref{observation1}(i), is equivalent to the $\{ {\rm claw},\overline{C_4}\}$-{\sc free Sandwich Problem}.
By Observation \ref{observation1}(iii), we may assume that $G_2$ has no universal vertex.
If $G_1$ is $\{ K_3,C_4\}$-free, then $G_1\in {\cal SW}_{\cal F}(G_1,G_2)$.
Hence, we may assume that $G_1$ contains an induced $K_3$ or $C_4$,
which implies that every graph in ${\cal SW}_{\cal F}(G_1,G_2)$ contains a triangle,
and, hence, in view of $\overline{\rm claw}$, is connected.

Suppose that ${\cal SW}_{\cal F}(G_1,G_2)$ contains some graph $G$ such that not all vertices of $G$ lie on triangles.
Let $T$ be the set of vertices of $G$ that lie on triangles, and let $R=V(G)\setminus T$,
in particular, $T$ and $R$ are both non-empty.
Since $G$ is $\overline{\rm claw}$-free, every vertex in $R$ has a neighbor in $T$.
If some vertex $u$ in $R$ has two neighbors $v$ and $w$ in $T$,
then, since $u$ does not lie on a triangle, $v$ and $w$ are not adjacent. 
Let $vxy$ be a triangle that contains $v$.
Since $G$ is $\overline{\rm claw}$-free, $w$ is adjacent to $x$ or $y$,
and $uvxwu$ or $uvywu$ is a $C_4$, which is a contradiction.
Hence, every vertex in $R$ has exactly one neighbor in $T$.
Let $v_1,\ldots,v_p$ be the vertices in $T$ that are the neighbor of some vertex in $R$.
Since $R$ is not empty, we have $p\geq 1$.
Let $u_i$ be a neighbor of $v_i$ in $R$ for $i\in [p]$.
Since $G$ is $\overline{\rm claw}$-free, every triangle of $G$ contains all vertices $v_1,\ldots,v_p$, which implies $p\leq 3$.
If $xy$ is an edge between two vertices in $R$, then $x$ and $y$ have different neighbors, say $v_i$ and $v_j$,
among $v_1,\ldots,v_p$. Since $v_i$ and $v_j$ both belong to every triangle, they are adjacent,
and the vertices $x$, $y$, $v_i$, and $v_j$ form a $C_4$, which is a contradiction. 
Hence, $R$ is independent.
If $p=1$, then $v_1$ is a universal vertex of $G$, and, hence, also of $G_2$,
which is a contradiction. Hence, $p\in \{ 2,3\}$.
If $p=3$, then $T=\{ v_1,v_2,v_3\}$, that is, $G$ contains exactly one triangle, and 
considering the $O\left(n(G_1)^3\right)$ choices for $v_1$, $v_2$, and $v_3$,
it is possible to check in polynomial time whether ${\cal SW}_{\cal F}(G_1,G_2)$ contains such a graph.
If $p=2$, then $V(G)\setminus \{ v_1,v_2\}$ is independent, which implies that 
considering the $O\left(n(G_1)^2\right)$ choices for $v_1$ and $v_2$,
it is possible to check in polynomial time whether ${\cal SW}_{\cal F}(G_1,G_2)$ contains such a graph.
Altogether, it follows that one can check in polynomial time 
whether ${\cal SW}_{\cal F}(G_1,G_2)$ contains some graph $G$ such that not all vertices of $G$ lie on triangles.
Hence, we may assume that all vertices of every graph in ${\cal SW}_{\cal F}(G_1,G_2)$ lie on triangles.

Suppose that $G$ is an edge-maximal graph in ${\cal SW}_{\cal F}(G_1,G_2)$.
If $G$ contains $\overline{K_1\cup C_4}$ as an induced subgraph, 
then the vertex of degree $4$ in $\overline{K_1\cup C_4}$ is universal in $G$, which is a contradiction.
Hence, $G$ is $\overline{K_1\cup C_4}$-free.
Our next goal is to show that $G$ contains the diamond as an induced subgraph.
Suppose, for a contradiction, that $G$ is diamond-free.
Since $G$ has no universal vertex, there is a triangle $uvw$ in $G$ as well as a vertex $x$ distinct from $u$, $v$, and $w$ 
such that $x$ in not adjacent to $u$.
Since $G$ is diamond-free, we may assume that $x$ is adjacent to $v$ but not to to $w$.
Since $x$ lies on some triangle, it has a neighbor $y$ outside of $\{ u,v,w\}$.
Since $G$ is $\left\{ \overline{\rm claw},C_4\right\}$-free, $y$ is adjacent to $v$.
Since $G$ is $\overline{K_1\cup C_4}$-free, $y$ is adjacent to $u$ or $w$, which yields a diamond in both cases.
Hence, $G$ contains a diamond.

Let $u$, $v$, $w$, and $x$ induce a diamond in $G$ such that $u$ and $x$ are not adjacent.
Clearly, every vertex of $G$ is adjacent to $v$ or $w$ or both.
Since $G_2$ has no universal vertex, 
the two sets $N_v=N_G(v)\setminus N_G[w]$ and $N_w=N_G(w)\setminus N_G[v]$ are both not empty.
Since $G$ is $\left\{ \overline{\rm claw},C_4\right\}$-free, $N_v\cup N_w$ is independent.
Let $R$ be the set of vertices in $N_G(v)\cap N_G(w)$ that have a neighbor in $N_v\cup N_w$, 
and let $S=(N_G(v)\cap N_G(w))\setminus R$.
By definition, $G$ contains no edge between $N_v\cup N_w$ and $S$.
Since every vertex of $G$ lies on a triangle, every vertex in $N_v\cup N_w$ has a neighbor in $R$,
in particular, $R$ is not empty.
Furthermore, since $G$ is $\overline{\rm claw}$-free, $N_v\cup N_w$ is completely joined to $R$.
If $R$ contains two non-adjacent vertices $x$ and $y$, then $x$, $y$, a vertex from $N_v$, and a vertex from $N_w$ form a $C_4$.
Hence, $R$ is a clique.
Since $G$ is $\overline{\rm claw}$-free, $S$ is independent.
%Since no vertex in $R$ is universal, every vertex in $R$ has a non-neighbor in $S$.
Since $G$ is $\overline{\rm claw}$-free, every vertex in $S$ has at most one non-neighbor in $R$.
Since adding an edge between $N_v$ and $w$ or between $N_w$ and $v$ 
does not create an induced subgraph $\overline{\rm claw}$ or $C_4$,
the edge-maximality of $G$ implies that $G$ contains all edges of $G_2$ between $\{ v,w\}$ and $V(G_1)\setminus \{ v,w\}$.
Similarly, since adding an edge between $R$ and $S$
does not create an induced subgraph $\overline{\rm claw}$ or $C_4$,
the edge-maximality of $G$ implies that $G$ contains all edges of $G_2$ between $R$ and $S$.
Altogether, it follows that there is an edge $vw$ of $G_2$, 
and a partition of the set $N_{G_2}(v)\cap N_{G_2}(w)$ into two sets $R$ and $S$ such that 
\begin{enumerate}[(i)]
\item the two sets $N_v=N_{G_2}(v)\setminus N_{G_2}[w]$ and $N_w=N_{G_2}(w)\setminus N_{G_2}[v]$ are non-empty,

$V(G_1)=\{ v,w\}\cup N_v\cup N_w\cup (N_{G_2}(v)\cap N_{G_2}(w))$, 

$N_v\cup N_w$ is independent in $G_1$, and,
\item $R$ is a clique in $G_2$,

$S$ is independent in $G_1$,

$G_2$ contains all possible edges between $R$ and $N_v\cup N_w$, 

$G_1$ contains no edge between $S$ and $N_v\cup N_w$, and,

in $G_2$, every vertex in $S$ has at most one non-neighbor in $R$.
\end{enumerate}
Conversely, if there is an edge $vw$ of $G_2$, 
and a partition of the set $N_{G_2}(v)\cap N_{G_2}(w)$ into two sets $R$ and $S$ such that (i) and (ii) are satisfied, 
then it is easy to see that ${\cal SW}_{\cal F}(G_1,G_2)$ is non-empty. 

Let $vw$ be an edge of $G_2$, and let $N_v$ and $N_w$ be as in (i).
Clearly, deciding whether (i) is satisfied can be done in polynomial time.
Furthermore, we now explain 
how to decide in polynomial time using $2${\sc Sat}
whether $N_{G_2}(v)\cap N_{G_2}(w)$ has a partition into two sets $R$ and $S$ that satisfies (ii).
Let $X=N_{G_2}(v)\cap N_{G_2}(w)$.
For every vertex $x$ in $X$, we introduce a boolean variable $x$,
which should be true if $x$ is in $R$, and false if $x$ is in $S$.
Now, we construct a $2${\sc Sat} formula $f$ as follows.
\begin{itemize}
\item For every two vertices $x$ and $y$ in $X$ that are non-adjacent in $G_2$, 
we add to $f$ the clause $\bar{x}\vee \bar{y}$,
reflecting that $R$ is a clique in $G_2$.
\item For every two vertices $x$ and $y$ in $X$ that are adjacent in $G_1$, 
we add to $f$ the clause $x\vee y$,
reflecting that $S$ is independent in $G_1$.
\item For every vertex $x$ in $X$ that is non-adjacent in $G_2$ to some vertex in $N_v\cup N_w$, 
we add to $f$ the clause $\bar{x}$,
reflecting that $G_2$ contains all possible edges between $R$ and $N_v\cup N_w$. 
\item For every vertex $x$ in $X$ that is adjacent in $G_1$ to some vertex in $N_v\cup N_w$, 
we add to $f$ the clause $x$,
reflecting that $G_1$ contains no edge between $S$ and $N_v\cup N_w$.
\item For every two vertices $x$ and $y$ in $X$, for which there is a third vertex $z$ in $X$ 
such that $x$ and $y$ are both non-adjacent in $G_2$ to $z$, 
we add to $f$ the clause $\bar{x}\vee \bar{y}$,
reflecting that, in $G_2$, every vertex in $S$ has at most one non-neighbor in $R$.
\end{itemize}
It is easy to see that $f$ is satisfiable if and only if 
$N_{G_2}(v)\cap N_{G_2}(w)$ has the desired partition.
Therefore, considering all $O\left(n(G_1)^2\right)$ edges $vw$ of $G_2$,
one can determine in polynomial time whether ${\cal SW}_{\cal F}(G_1,G_2)$ is non-empty.
$\Box$

\section{Some Hard Cases}

For every finite set ${\cal F}$ of graphs, 
the ${\cal F}$-{\sc free Sandwich Decision Problem} clearly belongs to NP.

Dantas, de Figueiredo, da Silva, and Teixeira \cite{dafisite} showed 
that the $\left\{ C_4\right\}$-{\sc free Sandwich Decision Problem} is NP-complete.
Considering the proof of the corresponding result (Theorem 1 in \cite{dafisite}),
it is easy to see that the very same proof yields the following result.

\begin{theorem}\label{theorem12}
The $\left\{ C_4,K_4\right\}$-{\sc free Sandwich Decision Problem} is NP-complete.
\end{theorem}
Our next two results rely on the hardness of deciding $3$-colorability.

\begin{theorem}\label{theorem11}
The $\left\{ {\rm paw},K_4\right\}$-{\sc free Sandwich Decision Problem} is NP-complete.
\end{theorem}
{\it Proof:} By Lemma \ref{lemma_olariu}, 
a connected graph $G$ that contains a triangle 
is $\left\{ {\rm paw},K_4\right\}$-free
if and only if 
it is $\left\{ \overline{P_3},K_4\right\}$-free
if and only if 
it is complete multipartite with at most three partite sets.
Furthermore, a graph is $3$-colorable if and only if it has a complete multipartite supergraph with at most three partite sets.
Therefore, a given connected graph $G$ that contains a triangle is $3$-colorable
if and only if ${\cal SW}_{\cal F}\left(G,K_{n(G)}\right)$ is non-empty for ${\cal F}=\left\{ {\rm paw},K_4\right\}$.
Since deciding $3$-colorability for such graphs is NP-complete,
the desired statement follows. $\Box$

\begin{theorem}\label{theorem16}
The $\left\{ {\rm paw},\overline{K_4}\right\}$-{\sc free Sandwich Decision Problem} is NP-complete.
\end{theorem}
{\it Proof:} Let $H$ be a graph.
Let $G_1$ arise from $H$ by adding three disjoint sets $X$, $Y$, and $Z$ each containing $R(3,4)$ new vertices,
and adding all edges between $X$ and $Y$, between $X$ and $Z$, and between $Y$ and $Z$.
Let $G_2$ arise from $G_1$ by adding all edges between $V(H)$ and $X\cup Y\cup Z$,
and by adding all edges of $\overline{H}$.
Note that $H$ is $3$-colorable if and only if $G_1$ is $3$-colorable.
Furthermore, for every $3$-coloring of $G_1$,
the three sets $X$, $Y$, and $Z$ are subsets of different color classes.
Hence, if $G_1$ is $3$-colorable, then $G_2$ contains a complete multipartite supergraph of $G_1$ with three partite sets.
We will show that $H$ is $3$-colorable if and only if ${\cal SW}_{\cal F}\left(G_1,G_2\right)$ is non-empty
for ${\cal F}=\left\{ \overline{\rm paw},K_4\right\}$.
By Observation \ref{observation1}(i), this implies the desired statement.

First, suppose that $H$ is $3$-colorable.
As observed above, $G_2$ contains a complete multipartite supergraph $G$ of $G_1$ with three partite sets.
Since $G$ is in ${\cal SW}_{\cal F}\left(G_1,G_2\right)$, the necessity follows.
For the proof of the sufficiency, suppose that $G$ is in ${\cal SW}_{\cal F}\left(G_1,G_2\right)$.
First, suppose that $\overline{G}$ is connected.
By Lemma \ref{lemma_olariu}, 
$\overline{G}$ is $\left\{ K_3,\overline{K_4}\right\}$-free or $\left\{ \overline{P_3},\overline{K_4}\right\}$-free.
In the first case, $n(G_1)\leq R(3,4)$, which is a contradiction.
In the second case, $\overline{G}$ is a complete multipartite graph with partite sets of order at most $3$,
which implies the contradiction $\Delta(G_1)\leq \Delta(G)\leq 2$.
Hence, $\overline{G}$ is disconnected.
Since $\overline{G}$ is $\overline{K_4}$-free, $\overline{G}$ has either two or three components.
First, suppose that $\overline{G}$ has two components.
Since $\overline{G}$ is $\left\{ {\rm paw},\overline{K_4}\right\}$-free,
one component, say $\overline{K}$, of $\overline{G}$ is a clique, and, by Lemma \ref{lemma_olariu},
the other component, say $\overline{K'}$, of $\overline{G}$ is 
$\left\{ K_3,\overline{K_3}\right\}$-free or $\left\{ \overline{P_3},\overline{K_3}\right\}$-free.
If $\overline{K'}$ is $\left\{ K_3,\overline{K_3}\right\}$-free, then $n(\overline{K'})\leq R(3,3)$.
Since $V\left(\overline{K}\right)$ is independent in $G$,
we may assume, by symmetry between $X$, $Y$, and $Z$, that $V\left(\overline{K'}\right)$ contains a vertex of $X$.
Since $G$ contains all edges between $V\left(\overline{K}\right)$ and $V\left(\overline{K'}\right)$,
this implies that the independent set $X$ is contained in $V\left(\overline{K}\right)$,
which is impossible because $|X|=R(3,4)>n(\overline{K'})$.
Hence, $\overline{K'}$ is $\left\{ \overline{P_3},\overline{K_3}\right\}$-free.
This implies that $\overline{K'}$ is a complete multipartite graph with partite sets of order at most $2$
and, hence, $G$, $G_1$, and $H$ are $3$-colorable.
Finally, suppose that $\overline{G}$ has three components.
It follows that each component of $\overline{G}$ is complete,
that is, $G$ is a complete multipartite graph with three partite sets.
Therefore, $G$, and, thus, also $G_1$ and $H$ are $3$-colorable.
$\Box$ 

\medskip

\noindent Our next two results rely on related results concerning $\left\{ \overline{C_4},C_3,C_5\right\}$-free graphs,
which are known as {\it chain graphs} or {\it difference graphs} (cf. Theorem 2.4.4 in \cite{mape}).
Clearly, chain graphs are bipartite.
While the next two proofs are based on essentially the same approach,
we argue from first principle for the first,
and rely on results about prime $\left\{ {\rm diamond},\overline{C_4}\right\}$-free graphs \cite{br} for the second.

\medskip

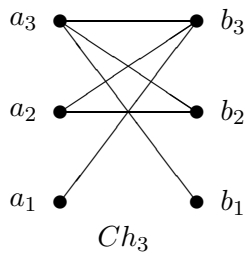
\begin{figure}[h]
\begin{center}
%TeXCAD Picture [4.pic]. Options:
%\grade{\on}
%\emlines{\off}
%\epic{\off}
%\beziermacro{\on}
%\reduce{\on}
%\snapping{\on}
%\pvinsert{% Your \input, \def, etc. here}
%\quality{8.000}
%\graddiff{0.005}
%\snapasp{1}
%\zoom{19.0276}
\unitlength 1.2mm % = 2.845pt
\linethickness{0.4pt}
\ifx\plotpoint\undefined\newsavebox{\plotpoint}\fi % GNUPLOT compatibility
\begin{picture}(23,25.5)(0,0)
\put(4,5){\circle*{1.5}}
\put(4,15){\circle*{1.5}}
\put(4,25){\circle*{1.5}}
\put(19,5){\circle*{1.5}}
\put(19,15){\circle*{1.5}}
\put(19,25){\circle*{1.5}}
\put(19,25){\line(-1,0){15}}
\put(4,25){\line(3,-2){15}}
\put(19,15){\line(-1,0){15}}
\put(4,15){\line(3,2){15}}
\put(19,25){\line(-3,-4){15}}
\put(19,5){\line(-3,4){15}}
\put(11,1){\makebox(0,0)[cc]{$Ch_3$}}
\put(0,25){\makebox(0,0)[cc]{$a_3$}}
\put(0,15){\makebox(0,0)[cc]{$a_2$}}
\put(0,5){\makebox(0,0)[cc]{$a_1$}}
\put(23,5){\makebox(0,0)[cc]{$b_1$}}
\put(23,15){\makebox(0,0)[cc]{$b_2$}}
\put(23,25){\makebox(0,0)[cc]{$b_3$}}
\end{picture}
\end{center}
\caption{The chain graph $Ch_3$.}\label{figchain}
\end{figure}

\begin{lemma}\label{lemma34}
Let ${\cal F}=\left\{ {\rm paw},\overline{C_4}\right\}$.
Let $G_1$ be the disjoint union of a graph $G$ and the graph $Ch_3$ in Figure \ref{figchain},
and let $G_2$ arise from $G_1$ by adding all edges between $V(G)$ and $V(Ch_3)$.
\begin{enumerate}[(i)]
\item If $G$ is $C_3$ or $C_5$, then ${\cal SW}_{\cal F}\left(G_1,G_2\right)$ is empty.
\item If $G$ is a chain graph, then ${\cal SW}_{\cal F}\left(G_1,G_2\right)$ contains a chain graph.
\end{enumerate}
\end{lemma}
{\it Proof:} (i) We only give details for the case that $G$ is a triangle $xyz$.
The case that $G$ is a $C_5$ can be settled similarly.

For a contradiction, suppose that ${\cal SW}_{\cal F}\left(G_1,G_2\right)$ contains a graph $H$.
Since $H$ is $\overline{C_4}$-free, considering the edges $a_3b_3$ and $xy$, 
we may assume, by symmetry, that $a_3$ and $x$ are adjacent.
Since $H$ is paw-free, considering $a_3$ and the triangle $xyz$,
we may assume, by symmetry, that $a_3$ is adjacent to $y$.
Since $H$ is paw-free, considering $b_3$ and the triangle $a_3xy$,
we may assume, by symmetry, that $b_3$ is adjacent to $y$.
Since $H$ is paw-free, considering any of the vertices $a_1$, $a_2$, and $b_2$ together with the triangle $a_3b_3y$,
we obtain that $a_1$, $a_2$, and $b_2$ are adjacent to $y$.
Now, $H[\{ a_1,a_2,b_2,y\}]$ is a paw, which is a contradiction.

\medskip

\noindent (ii) Let $G$ have the partite sets $A$ and $B$.
Let $G'$ arise from the disjoint union of $G$ and $Ch_3$
by adding all edges between $V(G)$ and $\{ a_1,a_2,a_3\}$.
Clearly, the sets $A'=A\cup \{ a_1,a_2,a_3\}$ and $B'=B\cup \{ b_1,b_2,b_3\}$ form a bipartition of $G'$.
Suppose that $G'$ contains an induced $\overline{C_4}$ with the two edges $ab$ and $a'b'$, where $a,a'\in A'$.
If $a\in \{ a_1,a_2,a_3\}$, then, in view of the edges between $V(G)$ and $\{ a_1,a_2,a_3\}$,
it follows that $b'\in \{ b_1,b_2,b_3\}$, which implies that $a'\in \{ a_1,a_2,a_3\}$,
and, hence, by symmetry, $b\in \{ b_1,b_2,b_3\}$.
Nevertheless, since $Ch_3$ is a chain graph, this is a contradiction.
If $a\in A$, then, in view of the structure of $G'$,
it follows that $b\in B$, which implies that $a'\in A$,
and, hence, by symmetry, $b'$.
Nevertheless, since $G$ is a chain graph, this is a contradiction.
Altogether, $G'$ is a bipartite $\overline{C_4}$-free graph,
that is, $G'$ is a chain graph.
By construction, $G'$ belongs to ${\cal SW}_{\cal F}\left(G_1,G_2\right)$.
$\Box$

\begin{theorem}\label{theorem34}
The $\left\{ {\rm paw},\overline{C_4}\right\}$-{\sc free Sandwich Decision Problem} is NP-complete.
\end{theorem}
{\it Proof:} Let ${\cal F}=\left\{ {\rm paw},\overline{C_4}\right\}$,
and let ${\Pi}$ be the set of all chain graphs.
In \cite{dafigoklma}
Dantas, Figueiredo, Golumbic, Klein, and Maffray
describe a polynomial reduction of an instance $f$ of an NP-complete variant of {\sc Satisfiability}
to an instance $(G_1,G_2)$ of the $\Pi$-{\sc Sandwich Decision Problem};
the decision variant of the $\Pi$-{\sc Sandwich Problem}.
Let $G_1'$ be the disjoint union of $G_1$ and the graph $Ch_3$, and, 
let $G_2'$ arise from the disjoint union of $G_2$ and the graph $Ch_3$ by adding all edges between $V(G_2)$ and $V(Ch_3)$.

If there is a chain graph $G$ with $G_1\subseteq G\subseteq  G_2$, 
then, by Lemma \ref{lemma34}(ii), ${\cal SW}_{\cal F}\left(G'_1,G'_2\right)$ is non-empty.
Conversely, if ${\cal SW}_{\cal F}\left(G'_1,G'_2\right)$ contains some graph $G'$,
then, by Lemma \ref{lemma34}(i),
the graph $G=G'-V(Ch_3)$ is $\left\{ \overline{C_4},C_3,C_5\right\}$-free,
that is, $G$ is a chain graph.
By construction, $G_1\subseteq G\subseteq  G_2$.
Altogether, we obtain a polynomial reduction of some NP-complete problem to the 
$\left\{ {\rm paw},\overline{C_4}\right\}$-{\sc free Sandwich Decision Problem},
which completes the proof. $\Box$

\begin{figure}[H]
\begin{center}
%TeXCAD Picture [6.pic]. Options:
%\grade{\on}
%\emlines{\off}
%\epic{\off}
%\beziermacro{\on}
%\reduce{\on}
%\snapping{\on}
%\pvinsert{% Your \input, \def, etc. here}
%\quality{8.000}
%\graddiff{0.005}
%\snapasp{1}
%\zoom{9.5137}
\unitlength 1mm % = 2.845pt
\linethickness{0.4pt}
\ifx\plotpoint\undefined\newsavebox{\plotpoint}\fi % GNUPLOT compatibility
\begin{picture}(30,54)(0,0)
\put(5,5){\circle*{1.5}}
\put(25,5){\circle*{1.5}}
\put(5,15){\circle*{1.5}}
\put(25,15){\circle*{1.5}}
\put(5,25){\circle*{1.5}}
\put(25,25){\circle*{1.5}}
\put(5,35){\circle*{1.5}}
\put(25,35){\circle*{1.5}}
\put(5,45){\circle*{1.5}}
\put(25,45){\circle*{1.5}}
\put(0,5){\makebox(0,0)[cc]{$b_1$}}
\put(30,5){\makebox(0,0)[cc]{$c_1$}}
\put(0,15){\makebox(0,0)[cc]{$b_2$}}
\put(30,15){\makebox(0,0)[cc]{$c_2$}}
\put(0,25){\makebox(0,0)[cc]{$b_3$}}
\put(30,25){\makebox(0,0)[cc]{$c_3$}}
\put(0,35){\makebox(0,0)[cc]{$b_4$}}
\put(30,35){\makebox(0,0)[cc]{$c_4$}}
\put(0,45){\makebox(0,0)[cc]{$b$}}
\put(30,45){\makebox(0,0)[cc]{$c$}}
\put(5,20){\oval(4,34)[]}
\put(25,45){\line(-5,-2){20}}
\put(25,20){\oval(4,34)[]}
\put(5,45){\line(5,-2){20}}
\put(15,50){\circle*{1.5}}
\put(15,50){\line(-2,-1){10}}
\put(5,45){\line(1,0){20}}
\put(25,45){\line(-2,1){10}}
\put(15,54){\makebox(0,0)[cc]{$a$}}
\put(25,35){\line(-1,0){20}}
\put(5,35){\line(2,-1){20}}
\put(25,25){\line(-1,0){20}}
\put(5,25){\line(2,1){20}}
\put(25,35){\line(-1,-1){20}}
\put(5,15){\line(2,1){20}}
\put(5,35){\line(1,-1){20}}
\put(5,35){\line(2,-3){20}}
\put(5,25){\line(2,-1){20}}
\put(5,5){\line(2,3){20}}
\put(15,0){\makebox(0,0)[cc]{$ECh_4$}}
\end{picture}
\end{center}
\caption{The graph $ECh_4$;
the vertex $b$ is adjacent to all vertices in the independent set $\{ c_1,c_2,c_3,c_4\}$, and, 
the vertex $c$ is adjacent to all vertices in the independent set $\{ b_1,b_2,b_3,b_4\}$.}\label{figech4}
\end{figure}
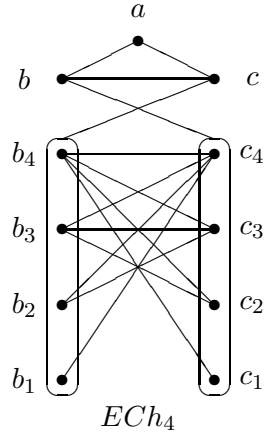

\begin{theorem}\label{theorem26}
The $\left\{ {\rm diamond},\overline{C_4}\right\}$-{\sc free Sandwich Decision Problem} is NP-complete.
\end{theorem}
{\it Proof:} Let ${\cal F}=\left\{ {\rm diamond},\overline{C_4}\right\}$,
and let ${\Pi}$ be the set of all chain graphs.
In \cite{dafigoklma}
Dantas, Figueiredo, Golumbic, Klein, and Maffray
describe a polynomial reduction of an instance $f$ of an NP-complete variant of {\sc Satisfiability}
to an instance $(G_1,G_2)$ of the $\Pi$-{\sc Sandwich Decision Problem},
where the edges of $G_1$ form a perfect matching.
We proceed similarly as in the proof of Theorem \ref{theorem34},
that is, we describe a polynomial reduction of $(G_1,G_2)$
to an instance $(G_1',G_2')$ of the ${\cal F}$-{\sc free Sandwich Decision Problem}.

Let $G_1'$ be the disjoint union of $G_1$ and the graph $ECh_4$ in Figure \ref{figech4}, and, 
let $G_2'$ arise from the disjoint union of $G_2$ and the graph $ECh_4$ 
by adding all edges between $V(G_2)$ and $V(ECh_4)\setminus \{ a\}$.

First, suppose that there is a chain graph $G$ with $G_1\subseteq G\subseteq G_2$.
Since $G_1$ has no isolated vertices, and $G$ is $\overline{C_4}$-free, 
it follows that $G$ is connected.
Let $B$ and $C$ be the partite sets of $G$.
Let $G'$ arise from the disjoint union of $G_1$ and the graph $ECh_4$
by adding all edges between $\{ b,b_1,b_2,b_3,b_4\}$ and $C$
as well as all edges between $\{ c\}$ and $B$.
Similarly as in the proof of Lemma \ref{lemma34}(ii), 
it follows that $G'\in {\cal SW}_{\cal F}(G'_1,G'_2)$.

Next, suppose that ${\cal SW}_{\cal F}(G'_1,G'_2)$ contains some graph $G'$.
Since $G'_1$ has no isolated vertices, and $G'$ is $\overline{C_4}$-free, 
it follows that $G'$ is connected.
In view of the edges of $Ech_4$, and, since the vertex $a$ has no neighbor in $V(G_1)$, 
also $\overline{G'}$ is connected.
Suppose that $U$ is a homogeneous set of $G'$.
Let $A$ be the set of vertices in $V(G)\setminus U$ that are adjacent to every vertex in $U$, and,
let $N$ be the set of vertices in $V(G)\setminus U$ that are adjacent to no vertex in $U$.
Since $G'$ and $\overline{G'}$ are connected, both sets $A$ and $N$ are non-empty. 
Since $G'$ is connected and diamond-free, the graph $G'[U]$ is $P_3$-free.
If $U$ is neither independent nor a clique, then, 
since $G'$ is diamond-free, it follows that $A$ contains only one vertex, and,
since $\overline{G'}$ is $\overline{C_4}$-free, it follows that $N$ is independent.
In this case, since $\overline{G}$ is connected, the unique vertex in $A$ is universal in $G'$,
which is a contradiction, because $G_2$ has no universal vertex.
If $U$ is a clique, then, since $G'$ is ${\cal F}$-free, it follows that $A$ is a clique and $N$ is independent.
Since $G'$ is diamond-free, every vertex in $N$ has exactly one neighbor in $A$.
In particular, it follows that the vertices of degree at least $2$ in $G'$ form a clique,
which is a contradiction in view of the two non-adjacent vertices $b_2$ and $c_2$, which have degree at least $2$ in $G_1$.
Hence, every homogeneous set of $G'$ is independent.
This easily implies that $\{ a\}$, $\{ b\}$, and $\{ c\}$ are maximal homogeneous sets of $G'$,
and, that all remaining vertices of $ECh_4$ belong to distinct maximal homogeneous sets of $G'$.
This implies that $Ech_4$ as an induced subgraph of the characteristic graph $G^*$ of $G'$;
in particular, the order of $G^*$ is at least $11$.
By a result of Brandst\"{a}dt \cite{br}, $G^*$ is 
\begin{enumerate}[(i)]
\item either a thin spider, that is, $V(G^*)$ can be partitioned into a clique $C$ and a stable set $S$,
and the edges between $C$ and $S$ form a matching that covers all of $S$ and all but at most one vertex of $C$,
\item or $G^*$ arises from the disjoint union of a triangle $a^*b^*c^*$ and a connected chain graph with partite sets $B^*$ and $C^*$
by adding all edges between $b^*$ and $C^*$ as well as all edges between $c^*$ and $B^*$.
\end{enumerate}
If $G^*$ is a thin spider, then, since the vertices that belong to maximal homogeneous sets represented by $S$ have independent neighborhoods, 
the vertices $a$, $b$, and $c$ of $ECh_4$ correspond to maximal homogeneous sets in $C$.
Since $b_4$ is non-adjacent to $c$, and, $c_4$ is non-adjacent to $b$,
the vertices $b_4$ and $c_4$ lie in $S$, which is a contradiction, since $S$ is independent.
Hence, $G^*$ is as in (ii).
Since $a$ is the only vertex of $G'$ whose removal yields a bipartite graph without creating a new vertex of degree $1$,
it follows that $\{ a\}=a^*$, and, by symmetry, $\{ b\}=b^*$ and $\{ c\}=c^*$.
Since $G^*-\{ a^*,b^*,c^*\}$ is a chain graph, and, every homogeneous set of $G'$ is independent,
it follows that $G'-\{ a,b,c\}$, and, hence also $G=G'-V(ECh_4)$ is a chain graph.
Since, by construction, $G_1\subseteq G\subseteq G_2$,
this completes the proof. $\Box$

\medskip

\noindent For our last two hardness results,
we prove the following auxiliary hardness result,
which might be of independent interest.

\begin{theorem}\label{theorem-co-matched}
Let $\Pi$ be the set of all bipartite graphs $G$ with a bipartition $A$ and $B$ such that 
every vertex in $A$ has at most one non-neighbor in $B$, and,
every vertex in $B$ has at most one non-neighbor in $A$.

The $\Pi$-{\sc Sandwich Decision Problem} is NP-complete.
\end{theorem}
{\it Proof:} The considered decision problem is clearly in NP.
In order to complete the proof, we describe a polynomial reduction 
of the well known NP-complete {\sc One-in-Three $3$Sat} (cf. [LO4] in \cite{gajo})
to the $\Pi$-{\sc Sandwich Decision Problem}.
Therefore, let $f$ be an instance of {\sc One-in-Three $3$Sat} 
consisting of the clauses $C_1,\ldots,C_m$
over the boolean variables $x_1,\ldots,x_n$.
We construct an instance $(G_1,G_2)$ of the $\Pi$-{\sc Sandwich Decision Problem} whose size is polynomially bounded in terms of $n$ and $m$
such that $f$ is a `yes'-instance of {\sc One-in-Three $3$Sat} 
if and only if $(G_1,G_2)$ is a `yes'-instance of the $\Pi$-{\sc Sandwich Decision Problem}.

\begin{figure}[H]
\begin{center}
%TeXCAD Picture [4.pic]. Options:
%\grade{\on}
%\emlines{\off}
%\epic{\off}
%\beziermacro{\on}
%\reduce{\on}
%\snapping{\on}
%\pvinsert{% Your \input, \def, etc. here}
%\quality{8.000}
%\graddiff{0.005}
%\snapasp{1}
%\zoom{9.5137}
\unitlength 0.9mm % = 2.845pt
\linethickness{0.4pt}
\ifx\plotpoint\undefined\newsavebox{\plotpoint}\fi % GNUPLOT compatibility
\begin{picture}(35,51)(0,0)
\put(5,5){\circle*{1.5}}
\put(30,5){\circle*{1.5}}
\put(5,20){\circle*{1.5}}
\put(30,20){\circle*{1.5}}
\put(5,35){\circle*{1.5}}
\put(30,35){\circle*{1.5}}
\put(5,50){\circle*{1.5}}
\put(30,50){\circle*{1.5}}
\put(0,50){\makebox(0,0)[cc]{$c^j$}}
\put(35,50){\makebox(0,0)[cc]{$d^j$}}
\put(-2,35){\makebox(0,0)[cc]{$p^j(x_1)$}}
\put(-2,20){\makebox(0,0)[cc]{$p^j(\bar{x}_2)$}}
\put(-2,5){\makebox(0,0)[cc]{$p^j(x_3)$}}
\put(35,35){\makebox(0,0)[cc]{$x_1^j$}}
\put(35,20){\makebox(0,0)[cc]{$\bar{x}_2^j$}}
\put(35,5){\makebox(0,0)[cc]{$x_3^j$}}
\put(5,5){\line(1,0){25}}
\put(5,20){\line(1,0){25}}
\put(5,35){\line(1,0){25}}
\put(5,50){\line(1,0){25}}
%\dashline{1}(5,50)(30,35)
\multiput(4.93,49.93)(.0537634,-.0322581){15}{\line(1,0){.0537634}}
\multiput(6.543,48.962)(.0537634,-.0322581){15}{\line(1,0){.0537634}}
\multiput(8.156,47.994)(.0537634,-.0322581){15}{\line(1,0){.0537634}}
\multiput(9.768,47.026)(.0537634,-.0322581){15}{\line(1,0){.0537634}}
\multiput(11.381,46.059)(.0537634,-.0322581){15}{\line(1,0){.0537634}}
\multiput(12.994,45.091)(.0537634,-.0322581){15}{\line(1,0){.0537634}}
\multiput(14.607,44.123)(.0537634,-.0322581){15}{\line(1,0){.0537634}}
\multiput(16.22,43.156)(.0537634,-.0322581){15}{\line(1,0){.0537634}}
\multiput(17.833,42.188)(.0537634,-.0322581){15}{\line(1,0){.0537634}}
\multiput(19.446,41.22)(.0537634,-.0322581){15}{\line(1,0){.0537634}}
\multiput(21.059,40.252)(.0537634,-.0322581){15}{\line(1,0){.0537634}}
\multiput(22.672,39.285)(.0537634,-.0322581){15}{\line(1,0){.0537634}}
\multiput(24.285,38.317)(.0537634,-.0322581){15}{\line(1,0){.0537634}}
\multiput(25.897,37.349)(.0537634,-.0322581){15}{\line(1,0){.0537634}}
\multiput(27.51,36.381)(.0537634,-.0322581){15}{\line(1,0){.0537634}}
\multiput(29.123,35.414)(.0537634,-.0322581){15}{\line(1,0){.0537634}}
%\end
%\dashline{1}(5,50)(30,20)
\multiput(4.93,49.93)(.0320924,-.0385109){19}{\line(0,-1){.0385109}}
\multiput(6.149,48.466)(.0320924,-.0385109){19}{\line(0,-1){.0385109}}
\multiput(7.369,47.003)(.0320924,-.0385109){19}{\line(0,-1){.0385109}}
\multiput(8.588,45.539)(.0320924,-.0385109){19}{\line(0,-1){.0385109}}
\multiput(9.808,44.076)(.0320924,-.0385109){19}{\line(0,-1){.0385109}}
\multiput(11.027,42.613)(.0320924,-.0385109){19}{\line(0,-1){.0385109}}
\multiput(12.247,41.149)(.0320924,-.0385109){19}{\line(0,-1){.0385109}}
\multiput(13.466,39.686)(.0320924,-.0385109){19}{\line(0,-1){.0385109}}
\multiput(14.686,38.222)(.0320924,-.0385109){19}{\line(0,-1){.0385109}}
\multiput(15.905,36.759)(.0320924,-.0385109){19}{\line(0,-1){.0385109}}
\multiput(17.125,35.296)(.0320924,-.0385109){19}{\line(0,-1){.0385109}}
\multiput(18.344,33.832)(.0320924,-.0385109){19}{\line(0,-1){.0385109}}
\multiput(19.564,32.369)(.0320924,-.0385109){19}{\line(0,-1){.0385109}}
\multiput(20.783,30.905)(.0320924,-.0385109){19}{\line(0,-1){.0385109}}
\multiput(22.003,29.442)(.0320924,-.0385109){19}{\line(0,-1){.0385109}}
\multiput(23.222,27.978)(.0320924,-.0385109){19}{\line(0,-1){.0385109}}
\multiput(24.442,26.515)(.0320924,-.0385109){19}{\line(0,-1){.0385109}}
\multiput(25.661,25.052)(.0320924,-.0385109){19}{\line(0,-1){.0385109}}
\multiput(26.881,23.588)(.0320924,-.0385109){19}{\line(0,-1){.0385109}}
\multiput(28.1,22.125)(.0320924,-.0385109){19}{\line(0,-1){.0385109}}
\multiput(29.32,20.661)(.0320924,-.0385109){19}{\line(0,-1){.0385109}}
%\end
%\dashline{1}(5,50)(30,5)
\multiput(4.93,49.93)(.0336927,-.0606469){14}{\line(0,-1){.0606469}}
\multiput(5.873,48.232)(.0336927,-.0606469){14}{\line(0,-1){.0606469}}
\multiput(6.817,46.533)(.0336927,-.0606469){14}{\line(0,-1){.0606469}}
\multiput(7.76,44.835)(.0336927,-.0606469){14}{\line(0,-1){.0606469}}
\multiput(8.703,43.137)(.0336927,-.0606469){14}{\line(0,-1){.0606469}}
\multiput(9.647,41.439)(.0336927,-.0606469){14}{\line(0,-1){.0606469}}
\multiput(10.59,39.741)(.0336927,-.0606469){14}{\line(0,-1){.0606469}}
\multiput(11.533,38.043)(.0336927,-.0606469){14}{\line(0,-1){.0606469}}
\multiput(12.477,36.345)(.0336927,-.0606469){14}{\line(0,-1){.0606469}}
\multiput(13.42,34.647)(.0336927,-.0606469){14}{\line(0,-1){.0606469}}
\multiput(14.364,32.949)(.0336927,-.0606469){14}{\line(0,-1){.0606469}}
\multiput(15.307,31.25)(.0336927,-.0606469){14}{\line(0,-1){.0606469}}
\multiput(16.25,29.552)(.0336927,-.0606469){14}{\line(0,-1){.0606469}}
\multiput(17.194,27.854)(.0336927,-.0606469){14}{\line(0,-1){.0606469}}
\multiput(18.137,26.156)(.0336927,-.0606469){14}{\line(0,-1){.0606469}}
\multiput(19.081,24.458)(.0336927,-.0606469){14}{\line(0,-1){.0606469}}
\multiput(20.024,22.76)(.0336927,-.0606469){14}{\line(0,-1){.0606469}}
\multiput(20.967,21.062)(.0336927,-.0606469){14}{\line(0,-1){.0606469}}
\multiput(21.911,19.364)(.0336927,-.0606469){14}{\line(0,-1){.0606469}}
\multiput(22.854,17.666)(.0336927,-.0606469){14}{\line(0,-1){.0606469}}
\multiput(23.798,15.967)(.0336927,-.0606469){14}{\line(0,-1){.0606469}}
\multiput(24.741,14.269)(.0336927,-.0606469){14}{\line(0,-1){.0606469}}
\multiput(25.684,12.571)(.0336927,-.0606469){14}{\line(0,-1){.0606469}}
\multiput(26.628,10.873)(.0336927,-.0606469){14}{\line(0,-1){.0606469}}
\multiput(27.571,9.175)(.0336927,-.0606469){14}{\line(0,-1){.0606469}}
\multiput(28.515,7.477)(.0336927,-.0606469){14}{\line(0,-1){.0606469}}
\multiput(29.458,5.779)(.0336927,-.0606469){14}{\line(0,-1){.0606469}}
%\end
%\dashline{1}(30,50)(5,35)
\multiput(29.93,49.93)(-.0537634,-.0322581){15}{\line(-1,0){.0537634}}
\multiput(28.317,48.962)(-.0537634,-.0322581){15}{\line(-1,0){.0537634}}
\multiput(26.704,47.994)(-.0537634,-.0322581){15}{\line(-1,0){.0537634}}
\multiput(25.091,47.026)(-.0537634,-.0322581){15}{\line(-1,0){.0537634}}
\multiput(23.478,46.059)(-.0537634,-.0322581){15}{\line(-1,0){.0537634}}
\multiput(21.865,45.091)(-.0537634,-.0322581){15}{\line(-1,0){.0537634}}
\multiput(20.252,44.123)(-.0537634,-.0322581){15}{\line(-1,0){.0537634}}
\multiput(18.639,43.156)(-.0537634,-.0322581){15}{\line(-1,0){.0537634}}
\multiput(17.026,42.188)(-.0537634,-.0322581){15}{\line(-1,0){.0537634}}
\multiput(15.414,41.22)(-.0537634,-.0322581){15}{\line(-1,0){.0537634}}
\multiput(13.801,40.252)(-.0537634,-.0322581){15}{\line(-1,0){.0537634}}
\multiput(12.188,39.285)(-.0537634,-.0322581){15}{\line(-1,0){.0537634}}
\multiput(10.575,38.317)(-.0537634,-.0322581){15}{\line(-1,0){.0537634}}
\multiput(8.962,37.349)(-.0537634,-.0322581){15}{\line(-1,0){.0537634}}
\multiput(7.349,36.381)(-.0537634,-.0322581){15}{\line(-1,0){.0537634}}
\multiput(5.736,35.414)(-.0537634,-.0322581){15}{\line(-1,0){.0537634}}
%\end
%\dashline{1}(30,50)(5,20)
\multiput(29.93,49.93)(-.0320924,-.0385109){19}{\line(0,-1){.0385109}}
\multiput(28.71,48.466)(-.0320924,-.0385109){19}{\line(0,-1){.0385109}}
\multiput(27.491,47.003)(-.0320924,-.0385109){19}{\line(0,-1){.0385109}}
\multiput(26.271,45.539)(-.0320924,-.0385109){19}{\line(0,-1){.0385109}}
\multiput(25.052,44.076)(-.0320924,-.0385109){19}{\line(0,-1){.0385109}}
\multiput(23.832,42.613)(-.0320924,-.0385109){19}{\line(0,-1){.0385109}}
\multiput(22.613,41.149)(-.0320924,-.0385109){19}{\line(0,-1){.0385109}}
\multiput(21.393,39.686)(-.0320924,-.0385109){19}{\line(0,-1){.0385109}}
\multiput(20.174,38.222)(-.0320924,-.0385109){19}{\line(0,-1){.0385109}}
\multiput(18.954,36.759)(-.0320924,-.0385109){19}{\line(0,-1){.0385109}}
\multiput(17.735,35.296)(-.0320924,-.0385109){19}{\line(0,-1){.0385109}}
\multiput(16.515,33.832)(-.0320924,-.0385109){19}{\line(0,-1){.0385109}}
\multiput(15.296,32.369)(-.0320924,-.0385109){19}{\line(0,-1){.0385109}}
\multiput(14.076,30.905)(-.0320924,-.0385109){19}{\line(0,-1){.0385109}}
\multiput(12.857,29.442)(-.0320924,-.0385109){19}{\line(0,-1){.0385109}}
\multiput(11.637,27.978)(-.0320924,-.0385109){19}{\line(0,-1){.0385109}}
\multiput(10.418,26.515)(-.0320924,-.0385109){19}{\line(0,-1){.0385109}}
\multiput(9.198,25.052)(-.0320924,-.0385109){19}{\line(0,-1){.0385109}}
\multiput(7.978,23.588)(-.0320924,-.0385109){19}{\line(0,-1){.0385109}}
\multiput(6.759,22.125)(-.0320924,-.0385109){19}{\line(0,-1){.0385109}}
\multiput(5.539,20.661)(-.0320924,-.0385109){19}{\line(0,-1){.0385109}}
%\end
%\dashline{1}(30,50)(5,5)
\multiput(29.93,49.93)(-.0336927,-.0606469){14}{\line(0,-1){.0606469}}
\multiput(28.986,48.232)(-.0336927,-.0606469){14}{\line(0,-1){.0606469}}
\multiput(28.043,46.533)(-.0336927,-.0606469){14}{\line(0,-1){.0606469}}
\multiput(27.1,44.835)(-.0336927,-.0606469){14}{\line(0,-1){.0606469}}
\multiput(26.156,43.137)(-.0336927,-.0606469){14}{\line(0,-1){.0606469}}
\multiput(25.213,41.439)(-.0336927,-.0606469){14}{\line(0,-1){.0606469}}
\multiput(24.269,39.741)(-.0336927,-.0606469){14}{\line(0,-1){.0606469}}
\multiput(23.326,38.043)(-.0336927,-.0606469){14}{\line(0,-1){.0606469}}
\multiput(22.383,36.345)(-.0336927,-.0606469){14}{\line(0,-1){.0606469}}
\multiput(21.439,34.647)(-.0336927,-.0606469){14}{\line(0,-1){.0606469}}
\multiput(20.496,32.949)(-.0336927,-.0606469){14}{\line(0,-1){.0606469}}
\multiput(19.552,31.25)(-.0336927,-.0606469){14}{\line(0,-1){.0606469}}
\multiput(18.609,29.552)(-.0336927,-.0606469){14}{\line(0,-1){.0606469}}
\multiput(17.666,27.854)(-.0336927,-.0606469){14}{\line(0,-1){.0606469}}
\multiput(16.722,26.156)(-.0336927,-.0606469){14}{\line(0,-1){.0606469}}
\multiput(15.779,24.458)(-.0336927,-.0606469){14}{\line(0,-1){.0606469}}
\multiput(14.835,22.76)(-.0336927,-.0606469){14}{\line(0,-1){.0606469}}
\multiput(13.892,21.062)(-.0336927,-.0606469){14}{\line(0,-1){.0606469}}
\multiput(12.949,19.364)(-.0336927,-.0606469){14}{\line(0,-1){.0606469}}
\multiput(12.005,17.666)(-.0336927,-.0606469){14}{\line(0,-1){.0606469}}
\multiput(11.062,15.967)(-.0336927,-.0606469){14}{\line(0,-1){.0606469}}
\multiput(10.118,14.269)(-.0336927,-.0606469){14}{\line(0,-1){.0606469}}
\multiput(9.175,12.571)(-.0336927,-.0606469){14}{\line(0,-1){.0606469}}
\multiput(8.232,10.873)(-.0336927,-.0606469){14}{\line(0,-1){.0606469}}
\multiput(7.288,9.175)(-.0336927,-.0606469){14}{\line(0,-1){.0606469}}
\multiput(6.345,7.477)(-.0336927,-.0606469){14}{\line(0,-1){.0606469}}
\multiput(5.401,5.779)(-.0336927,-.0606469){14}{\line(0,-1){.0606469}}
%\end
%\dashline{1}(30,35)(5,20)
\multiput(29.93,34.93)(-.0537634,-.0322581){15}{\line(-1,0){.0537634}}
\multiput(28.317,33.962)(-.0537634,-.0322581){15}{\line(-1,0){.0537634}}
\multiput(26.704,32.994)(-.0537634,-.0322581){15}{\line(-1,0){.0537634}}
\multiput(25.091,32.026)(-.0537634,-.0322581){15}{\line(-1,0){.0537634}}
\multiput(23.478,31.059)(-.0537634,-.0322581){15}{\line(-1,0){.0537634}}
\multiput(21.865,30.091)(-.0537634,-.0322581){15}{\line(-1,0){.0537634}}
\multiput(20.252,29.123)(-.0537634,-.0322581){15}{\line(-1,0){.0537634}}
\multiput(18.639,28.156)(-.0537634,-.0322581){15}{\line(-1,0){.0537634}}
\multiput(17.026,27.188)(-.0537634,-.0322581){15}{\line(-1,0){.0537634}}
\multiput(15.414,26.22)(-.0537634,-.0322581){15}{\line(-1,0){.0537634}}
\multiput(13.801,25.252)(-.0537634,-.0322581){15}{\line(-1,0){.0537634}}
\multiput(12.188,24.285)(-.0537634,-.0322581){15}{\line(-1,0){.0537634}}
\multiput(10.575,23.317)(-.0537634,-.0322581){15}{\line(-1,0){.0537634}}
\multiput(8.962,22.349)(-.0537634,-.0322581){15}{\line(-1,0){.0537634}}
\multiput(7.349,21.381)(-.0537634,-.0322581){15}{\line(-1,0){.0537634}}
\multiput(5.736,20.414)(-.0537634,-.0322581){15}{\line(-1,0){.0537634}}
%\end
%\dashline{1}(5,20)(30,5)
\multiput(4.93,19.93)(.0537634,-.0322581){15}{\line(1,0){.0537634}}
\multiput(6.543,18.962)(.0537634,-.0322581){15}{\line(1,0){.0537634}}
\multiput(8.156,17.994)(.0537634,-.0322581){15}{\line(1,0){.0537634}}
\multiput(9.768,17.026)(.0537634,-.0322581){15}{\line(1,0){.0537634}}
\multiput(11.381,16.059)(.0537634,-.0322581){15}{\line(1,0){.0537634}}
\multiput(12.994,15.091)(.0537634,-.0322581){15}{\line(1,0){.0537634}}
\multiput(14.607,14.123)(.0537634,-.0322581){15}{\line(1,0){.0537634}}
\multiput(16.22,13.156)(.0537634,-.0322581){15}{\line(1,0){.0537634}}
\multiput(17.833,12.188)(.0537634,-.0322581){15}{\line(1,0){.0537634}}
\multiput(19.446,11.22)(.0537634,-.0322581){15}{\line(1,0){.0537634}}
\multiput(21.059,10.252)(.0537634,-.0322581){15}{\line(1,0){.0537634}}
\multiput(22.672,9.285)(.0537634,-.0322581){15}{\line(1,0){.0537634}}
\multiput(24.285,8.317)(.0537634,-.0322581){15}{\line(1,0){.0537634}}
\multiput(25.897,7.349)(.0537634,-.0322581){15}{\line(1,0){.0537634}}
\multiput(27.51,6.381)(.0537634,-.0322581){15}{\line(1,0){.0537634}}
\multiput(29.123,5.414)(.0537634,-.0322581){15}{\line(1,0){.0537634}}
%\end
%\dashline{1}(30,5)(5,35)
\multiput(29.93,4.93)(-.0320924,.0385109){19}{\line(0,1){.0385109}}
\multiput(28.71,6.393)(-.0320924,.0385109){19}{\line(0,1){.0385109}}
\multiput(27.491,7.857)(-.0320924,.0385109){19}{\line(0,1){.0385109}}
\multiput(26.271,9.32)(-.0320924,.0385109){19}{\line(0,1){.0385109}}
\multiput(25.052,10.783)(-.0320924,.0385109){19}{\line(0,1){.0385109}}
\multiput(23.832,12.247)(-.0320924,.0385109){19}{\line(0,1){.0385109}}
\multiput(22.613,13.71)(-.0320924,.0385109){19}{\line(0,1){.0385109}}
\multiput(21.393,15.174)(-.0320924,.0385109){19}{\line(0,1){.0385109}}
\multiput(20.174,16.637)(-.0320924,.0385109){19}{\line(0,1){.0385109}}
\multiput(18.954,18.1)(-.0320924,.0385109){19}{\line(0,1){.0385109}}
\multiput(17.735,19.564)(-.0320924,.0385109){19}{\line(0,1){.0385109}}
\multiput(16.515,21.027)(-.0320924,.0385109){19}{\line(0,1){.0385109}}
\multiput(15.296,22.491)(-.0320924,.0385109){19}{\line(0,1){.0385109}}
\multiput(14.076,23.954)(-.0320924,.0385109){19}{\line(0,1){.0385109}}
\multiput(12.857,25.418)(-.0320924,.0385109){19}{\line(0,1){.0385109}}
\multiput(11.637,26.881)(-.0320924,.0385109){19}{\line(0,1){.0385109}}
\multiput(10.418,28.344)(-.0320924,.0385109){19}{\line(0,1){.0385109}}
\multiput(9.198,29.808)(-.0320924,.0385109){19}{\line(0,1){.0385109}}
\multiput(7.978,31.271)(-.0320924,.0385109){19}{\line(0,1){.0385109}}
\multiput(6.759,32.735)(-.0320924,.0385109){19}{\line(0,1){.0385109}}
\multiput(5.539,34.198)(-.0320924,.0385109){19}{\line(0,1){.0385109}}
%\end
%\dashline{1}(5,35)(30,20)
\multiput(4.93,34.93)(.0537634,-.0322581){15}{\line(1,0){.0537634}}
\multiput(6.543,33.962)(.0537634,-.0322581){15}{\line(1,0){.0537634}}
\multiput(8.156,32.994)(.0537634,-.0322581){15}{\line(1,0){.0537634}}
\multiput(9.768,32.026)(.0537634,-.0322581){15}{\line(1,0){.0537634}}
\multiput(11.381,31.059)(.0537634,-.0322581){15}{\line(1,0){.0537634}}
\multiput(12.994,30.091)(.0537634,-.0322581){15}{\line(1,0){.0537634}}
\multiput(14.607,29.123)(.0537634,-.0322581){15}{\line(1,0){.0537634}}
\multiput(16.22,28.156)(.0537634,-.0322581){15}{\line(1,0){.0537634}}
\multiput(17.833,27.188)(.0537634,-.0322581){15}{\line(1,0){.0537634}}
\multiput(19.446,26.22)(.0537634,-.0322581){15}{\line(1,0){.0537634}}
\multiput(21.059,25.252)(.0537634,-.0322581){15}{\line(1,0){.0537634}}
\multiput(22.672,24.285)(.0537634,-.0322581){15}{\line(1,0){.0537634}}
\multiput(24.285,23.317)(.0537634,-.0322581){15}{\line(1,0){.0537634}}
\multiput(25.897,22.349)(.0537634,-.0322581){15}{\line(1,0){.0537634}}
\multiput(27.51,21.381)(.0537634,-.0322581){15}{\line(1,0){.0537634}}
\multiput(29.123,20.414)(.0537634,-.0322581){15}{\line(1,0){.0537634}}
%\end
%\dashline{1}(30,20)(5,5)
\multiput(29.93,19.93)(-.0537634,-.0322581){15}{\line(-1,0){.0537634}}
\multiput(28.317,18.962)(-.0537634,-.0322581){15}{\line(-1,0){.0537634}}
\multiput(26.704,17.994)(-.0537634,-.0322581){15}{\line(-1,0){.0537634}}
\multiput(25.091,17.026)(-.0537634,-.0322581){15}{\line(-1,0){.0537634}}
\multiput(23.478,16.059)(-.0537634,-.0322581){15}{\line(-1,0){.0537634}}
\multiput(21.865,15.091)(-.0537634,-.0322581){15}{\line(-1,0){.0537634}}
\multiput(20.252,14.123)(-.0537634,-.0322581){15}{\line(-1,0){.0537634}}
\multiput(18.639,13.156)(-.0537634,-.0322581){15}{\line(-1,0){.0537634}}
\multiput(17.026,12.188)(-.0537634,-.0322581){15}{\line(-1,0){.0537634}}
\multiput(15.414,11.22)(-.0537634,-.0322581){15}{\line(-1,0){.0537634}}
\multiput(13.801,10.252)(-.0537634,-.0322581){15}{\line(-1,0){.0537634}}
\multiput(12.188,9.285)(-.0537634,-.0322581){15}{\line(-1,0){.0537634}}
\multiput(10.575,8.317)(-.0537634,-.0322581){15}{\line(-1,0){.0537634}}
\multiput(8.962,7.349)(-.0537634,-.0322581){15}{\line(-1,0){.0537634}}
\multiput(7.349,6.381)(-.0537634,-.0322581){15}{\line(-1,0){.0537634}}
\multiput(5.736,5.414)(-.0537634,-.0322581){15}{\line(-1,0){.0537634}}
%\end
%\dashline{1}(5,5)(30,35)
\multiput(4.93,4.93)(.0320924,.0385109){19}{\line(0,1){.0385109}}
\multiput(6.149,6.393)(.0320924,.0385109){19}{\line(0,1){.0385109}}
\multiput(7.369,7.857)(.0320924,.0385109){19}{\line(0,1){.0385109}}
\multiput(8.588,9.32)(.0320924,.0385109){19}{\line(0,1){.0385109}}
\multiput(9.808,10.783)(.0320924,.0385109){19}{\line(0,1){.0385109}}
\multiput(11.027,12.247)(.0320924,.0385109){19}{\line(0,1){.0385109}}
\multiput(12.247,13.71)(.0320924,.0385109){19}{\line(0,1){.0385109}}
\multiput(13.466,15.174)(.0320924,.0385109){19}{\line(0,1){.0385109}}
\multiput(14.686,16.637)(.0320924,.0385109){19}{\line(0,1){.0385109}}
\multiput(15.905,18.1)(.0320924,.0385109){19}{\line(0,1){.0385109}}
\multiput(17.125,19.564)(.0320924,.0385109){19}{\line(0,1){.0385109}}
\multiput(18.344,21.027)(.0320924,.0385109){19}{\line(0,1){.0385109}}
\multiput(19.564,22.491)(.0320924,.0385109){19}{\line(0,1){.0385109}}
\multiput(20.783,23.954)(.0320924,.0385109){19}{\line(0,1){.0385109}}
\multiput(22.003,25.418)(.0320924,.0385109){19}{\line(0,1){.0385109}}
\multiput(23.222,26.881)(.0320924,.0385109){19}{\line(0,1){.0385109}}
\multiput(24.442,28.344)(.0320924,.0385109){19}{\line(0,1){.0385109}}
\multiput(25.661,29.808)(.0320924,.0385109){19}{\line(0,1){.0385109}}
\multiput(26.881,31.271)(.0320924,.0385109){19}{\line(0,1){.0385109}}
\multiput(28.1,32.735)(.0320924,.0385109){19}{\line(0,1){.0385109}}
\multiput(29.32,34.198)(.0320924,.0385109){19}{\line(0,1){.0385109}}
%\end
\end{picture}
\end{center}
\caption{The vertices added for the clause $C_j=x_1\vee \bar{x}_2\vee x_3$. 
The dashed lines are the edges in $E_j$.
}\label{fig-co-matched}
\end{figure}
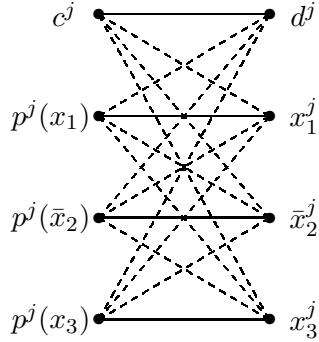

Starting with the empty graph, we construct $G_1$ as follows.
\begin{itemize}
\item For every clause $C_j$ with literals $u$, $v$, and $w$, 
add the eight vertices $c^j$, $d^j$, $u^j$, $v^j$, $w^j$, $p^j(u)$, $p^j(v)$, and $p^j(w)$,
add the four edges $c^jd^j$, $u^jp(u)^j$, $v^jp(v)^j$, and $w^jp(w)^j$,
and let 
\begin{eqnarray*}
E_j & = & \{ c^ju^j,c^jv^j,c^jw^j\}\cup \{ d^jp^j(u),d^jp^j(v),d^jp^j(w)\}\\
&& \cup \{u^jp^j(v),u^jp^j(w),v^jp^j(u),v^jp^j(w),w^jp^j(u),w^jp^j(v)\}.
\end{eqnarray*}
See Figure \ref{fig-co-matched} for an example.
\item For every $i,j\in [m]$, and $k\in [n]$, add the edge $c^id^j$, and, if the corresponding vertices exists, the edge $x_k^i\bar{x}_k^j$.
\end{itemize}
Let $G_2$ arise from $G_1$ by adding all edges in $E\left(\overline{G_1}\right)\setminus\bigcup\limits_{j\in [m]}E_j$.
Clearly, the size of $(G_1,G_2)$ is polynomially bounded in terms of $n$ and $m$.

First, suppose that there is a truth assignment such that every clause of $f$ contains exactly one true literal.
Let
\begin{eqnarray*}
T & = & 
\bigcup_{j\in [m]}
\Big\{ d^j\Big\}
\cup\bigcup_{i\in [n]:x_i\,\,is\,\,a\,\,true\,\,literal\,\,in\,\,C_j}\Big\{ x_i^j\Big\}
\cup\bigcup_{i\in [n]:x_i\,\,is\,\,a\,\,false\,\,literal\,\,in\,\,C_j}\Big\{ p(x_i)^j\Big\},\\
F & = & 
\bigcup_{j\in [m]}
\Big\{ c^j\Big\}
\cup\bigcup_{i\in [n]:x_i\,\,is\,\,a\,\,true\,\,literal\,\,in\,\,C_j}\Big\{ p(x_i)^j\Big\}
\cup\bigcup_{i\in [n]:x_i\,\,is\,\,a\,\,false\,\,literal\,\,in\,\,C_j}\Big\{ x_i^j\Big\},\\
V(G) &=& V(G_1)\mbox{, and}\\
E(G) & = & \{ e\in E(G_2):|e\cap T|=1\}.  
\end{eqnarray*}
Clearly, $G_1\subseteq G\subseteq G_2$.
Furthermore,
if the clause $C_j$ contains the three literals $u$, $v$, and $w$, and $u$ is the true literal in $C_j$, then, within the graph $G$,
\begin{itemize}
\item $u^j$ is the only non-neighbor of $c^j$ in $T$,
$c^j$ is the only non-neighbor of $u^j$ in $F$,
\item $p(u)^j$ is the only non-neighbor of $d^j$ in $F$,
$d^j$ is the only non-neighbor of $p(u)^j$ in $T$,
\item $p(w)^j$ is the only non-neighbor of $v^j$ in $T$, 
$v^j$ is the only non-neighbor of $p(w)^j$ in $F$, 
\item $p(v)^j$ is the only non-neighbor of $w^j$ in $T$, and
$w^j$ is the only non-neighbor of $p(v)^j$ in $F$.
\end{itemize}
Altogether, it follows that $G$ solves the $\Pi$-{\sc Sandwich Problem}.

Now, suppose that $G$ solves the $\Pi$-{\sc Sandwich Problem}.
Let $T$ and $F$ denote a bipartition of $G$.
In view of the edges of $G_1$, 
we may assume that $T$ contains the vertices $d^1,\ldots,d^m$,
and, that $F$ contains the vertices $c^1,\ldots,c^m$.
Let $j\in [m]$, and let the clause $C_j$ contain the literals $u$, $v$, and $w$.
In view of the edges in $E_j$ incident with $c^j$ and $d^j$, respectively, it follows 
that $T$ contains at most one vertex from $\{ u,v,w\}$, and,
that $F$ contains at most one vertex from $\{ p(u),p(v),p(w)\}$.
Now, the edges in $E_j$ between $\{ u,v,w\}$ and $\{ p(u),p(v),p(w)\}$ imply 
that $T$ contains exactly one vertex from $\{ u,v,w\}$.
Hence, in view of the edges of the form $x_k^i\bar{x}_k^j$, 
setting the variables $x_i$ that correspond to a vertex $x_i^j$ in $T$ to true
yields a consistent truth assignment for which each clause of $f$ contains exactly one true literal.
$\Box$

\begin{theorem}\label{theorem25}
The $\left\{ {\rm paw},\overline{\rm diamond}\right\}$-{\sc free Sandwich Decision Problem} is NP-complete.
\end{theorem}
{\it Proof:} Let ${\cal F}=\left\{ {\rm paw},\overline{\rm diamond}\right\}$,
and, let $\Pi$ be as in Theorem \ref{theorem-co-matched}.
We describe a polynomial reduction of an instance $(G_1,G_2)$ of the NP-complete $\Pi$-{\sc Sandwich Decision Problem}
to an instance $(G_1',G_2')$ of the ${\cal F}$-{\sc free Sandwich Decision Problem}.
In view of the proof of Theorem \ref{theorem-co-matched},
we may assume that $G_1$ has order at least $8$, and contains no isolated vertex.

Let $P:a_1b_1a_2b_2$ be an induced $P_4$.
Let $G_1'$ be the disjoint union of $G_1$ and $P$, and, 
let $G_2'$ arise from the disjoint union of $G_2$ and $P$ 
by adding all edges between $V(G_2)$ and $V(P)$.

First, suppose that $G$ solves the $\Pi$-{\sc Sandwich Problem},
and, that the sets $A$ and $B$ form a suitable bipartition of $G$.
If $G'$ arises from the disjoint union of $G$ and $P$
by adding all edges between $\{ a_1,a_2\}$ and $B$
as well as all edges between $\{ b_1,b_2\}$ and $A$,
then $G'\in {\cal SW}_{\cal F}\left(G'_1,G'_2\right)$.
Conversely, suppose that ${\cal SW}_{\cal F}\left(G'_1,G'_2\right)$ contains a graph $G'$.
In view of $P$, some component of $G'$ contains an induced $\overline{P_3}$.
Since $G'$ is $\overline{\rm diamond}$-free, this implies that $G$ is connected.
By Lemma \ref{lemma_olariu}, $G'$ is $\left\{ K_3,\overline{\rm diamond}\right\}$-free.
Suppose that $G'$ is not bipartite.
Let $C:u_1\ldots u_\ell$ be a shortest odd cycle in $G'$.
Since $G'$ is triangle-free, $\ell$ is at least $5$.
Since $G'$ is $\overline{\rm diamond}$-free, $\ell$ is at most $5$,
that is, $\ell$ is $5$.
Since $G'$ has order more than $5$, there is some vertex $v$ in $V(G')\setminus V(C)$.
Since $G'$ is triangle-free, we may assume, by symmetry, that $N_G(v)\cap V(C)$ is contained in $\{ u_1,u_3\}$.
Now, $G'[\{ u_2,u_4,u_5,v\}]$ is a $\overline{\rm diamond}$, which is a contradiction.
Hence, $G'$ is bipartite. 
Let the sets $A'$ and $B'$ form a bipartition of $G'$ with $a_1,a_2\in A'$ and $b_1,b_2\in B'$.
Let $G=G'-V(P)$.
Let $A=A'\setminus \{ a_1,a_2\}$ and $B=B'\setminus \{ b_1,b_2\}$.
Suppose that some vertex $a$ in $A$ is non-adjacent to two vertices $b$ and $b'$ in $B$.
Since $G_1$ has no isolated vertex, $a$ has a neighbor $b''$ in $B$, and, 
$G[\{ a,b,b',b''\}]$ is a $\overline{\rm diamond}$,
which is a contradiction.
By symmetry, it follows that $G$ solves the $\Pi$-{\sc Sandwich Problem}
$\Box$

\begin{theorem}\label{theorem27}
The $\left\{ {\rm diamond},\overline{\rm diamond}\right\}$-{\sc free Sandwich Decision Problem} is NP-complete.
\end{theorem}
{\it Proof:} Let ${\cal F}=\left\{ {\rm diamond},\overline{\rm diamond}\right\}$,
and, let $\Pi$ be as in Theorem \ref{theorem-co-matched}.
We describe a polynomial reduction of an instance $(G_1,G_2)$ of the NP-complete $\Pi$-{\sc Sandwich Decision Problem}
to an instance $(G_1',G_2')$ of the ${\cal F}$-{\sc free Sandwich Decision Problem}.
In view of the proof of Theorem \ref{theorem-co-matched},	
we may assume that $G_1$ has order at least $8$, and contains no isolated vertex.

Let $P'$ be the graph with vertices $a_1$, $a_2$, $b_1$, $b_1'$, and $b_2$,
and edges $a_1b_1$, $a_1b_1'$, $b_1,a_2$, and $a_2b_2$.
Let $G_1'$ be the disjoint union of $G_1$ and $P'$, and, 
let $G_2'$ arise from the disjoint union of $G_2$ and $P'$ 
by adding all edges between $V(G_2)$ and $V(P')$.

First, suppose that $G$ solves the $\Pi$-{\sc Sandwich Problem},
and, that the sets $A$ and $B$ form a suitable bipartition of $G$.
If $G'$ arises from the disjoint union of $G$ and $P'$
by adding all edges between $\{ a_1,a_2\}$ and $B$
as well as all edges between $\{ b_1,b_1',b_2\}$ and $A$,
then $G'\in {\cal SW}_{\cal F}\left(G'_1,G'_2\right)$.
Conversely, suppose that ${\cal SW}_{\cal F}\left(G'_1,G'_2\right)$ contains a graph $G'$.
Since $P$ contains an induced $P_4$, 
some component of $G'$ as well as some component of $\overline{G'}$ contains an induced $\overline{P_3}$.
Since $G'$ and $\overline{G'}$ are $\overline{\rm diamond}$-free, this implies that $G$ and $\overline{G'}$ are connected.
If $G'$ is prime, then, since $\overline{G'}$ contains the triangle $b_1b_1'b_2$,
a result of Brandst\"{a}dt and Mahfud \cite{brma} implies 
that $G'$ is a bipartite graph with partite sets $A'$ and $B'$ such that 
every vertex in $A'$ has at most one non-neighbor in $B'$, and,
every vertex in $B'$ has at most one non-neighbor in $A'$.
Now, $G=G'-V(P')$ solves the $\Pi$-{\sc Sandwich Problem}.
Hence, we may assume that $G'$ contains a homogeneous set $U$ of vertices,
that is, $2\leq |U|\leq n(G')-1$, and $V(G')=U\cup A\cup B$,
where $A$ is the set of vertices in $V(G')\setminus U$ that are adjacent to every vertex in $U$, and,
$B$ is the set of vertices in $V(G')\setminus U$ that are adjacent to no vertex in $U$.
Suppose that $U$ contains two adjacent vertices.
Since $G'$ is diamond-free, $A$ is a clique.
Since $G'$ is $\overline{\rm diamond}$-free, $B$ is a clique.
Since $G'$ and $\overline{G'}$ are connected,
we have $|A|\geq 2$ or $|B|\geq 2$.
Since $G'$ is ${\cal F}$-free, $U$ is a clique.
Nevertheless, this implies that $\overline{G'}$ is bipartite,
which is impossible in view of the triangle $b_1b_1'b_2$ in $\overline{G'}$.
Hence, we may assume that $U$ contains two non-adjacent vertices.
Arguing similarly as above, this implies that $A$, $B$, and $U$ are independent,
that is, $G'$ is bipartite with bipartition $A$ and $U\cup B$.
If some vertex $a$ in $A$ has two non-neighbors $b$ and $b'$ in $B$,
then $a$, $b$, $b'$, and a vertex from $U$ induce a $\overline{\rm diamond}$,
which is a contradiction.
If some vertex $b$ in $B$ has two non-neighbors $a$ and $a'$ in $A$,
then, since $G$ is connected, $b$ has a neighbor $a''$ in $A$, and,
$b$, $a$, $a'$, and $a''$ induce a $\overline{\rm diamond}$,
which is a contradiction.
This completes the proof. 
$\Box$

\section{Conclusion}

Figure \ref{fig1} shows eight open cases, and, since the hardness results were typically slightly harder to obtain,
we tend to believe that most of the corresponding problems are hard.

\end{document}